\documentclass{amsart}
\usepackage{amsthm}
\newtheorem{lemma}{Lemma}

\newtheorem{theorem}{Theorem}
\newtheorem{corollary}{Corollary}
\newtheorem{definition}{Definition}

\newtheorem{remark}{Remark}
\newtheorem{example}{Example}
\newtheorem{examples}{Examples}
\usepackage{changebar}

\newtheorem{fact}{Fact}

\newcommand{\nth}[1]{$#1 {\rm - th }$}

\newcommand{\rad}{{\rm rad } \ }
\newcommand{\lcm}{\mbox{\ lcm\ }}
\newcommand{\Z}{\mathbb{Z}} \newcommand{\ZM}[1]{\Z /( #1 \cdot \Z)}
\newcommand{\ZMs}[1]{\left(\Z / #1 \cdot \Z\right)^{\ast}}
\newcommand{\ZMb}[1]{\left(\Z / #1 \cdot \Z\right)^{\top}}
\newcommand{\ZMd}[1]{\left(\Z / #1 \cdot \Z\right)^{\widehat{}} \ }
\newcommand{\ord}{{\rm ord}} \newcommand{\rem}{{\rm rem }}
 \newcommand{\disc}{{\rm disc}}
\newcommand{\Tr}{\mbox{\bf Tr}}

\newcommand{\rg}[1]{\mbox{\bf #1}}
\newcommand{\eu}[1]{\mathfrak{#1}}
\newcommand{\id}[1]{\mathcal{#1}}
\newcommand{\Gal}{\mbox{Gal}}
\newcommand{\Aut}{\mbox{Aut}}

\newcommand{\rf}[1]{(\ref{#1})}
\newcommand{\smod}{\mbox{\scriptsize \ mod\ }}

\newcommand{\lchooses}[2]{\left( \frac{#1}{#2 } \right)}

\newcommand{\F}{\mathbb{F}}
\newcommand{\K}{\mathbb{K}}
\newcommand{\LK}{\mathbb{L}}
\newcommand{\KL}{\mathbb{L}}

\newcommand{\Q}{\mathbb{Q}}

\newcommand{\C}{\mathbb{C}}
\newcommand{\N}{\mathbb{N}}

\newcommand{\wh}{\widehat}
\def\NI{\noindent}

\begin{document}

\title{Cyclotomy Primality Proofs and their Certificates} 
\author{Preda Mih\u{a}ilescu}
\address{Mathematisches Institut der Universit\"at
G\"ottingen}
\email{preda@uni-amth.gwdg.de} 

\date{Version 2.0 \today}
\thanks{The author is sponored by the Volkswagen Foundation}
\bigskip
{\obeylines \em 
\vspace*{0.8cm}
\noindent\hspace*{5cm}Elle est \`a toi cette chanson\newline
\vspace*{-0.4cm}
\noindent\hspace*{5cm}Toi l'professeur qui sans fa\c{c}on, \newline
\vspace*{-0.4cm}
\noindent\hspace*{5cm}As ouvert ma petite th\`ese \newline
\vspace*{-0.4cm}
\noindent\hspace*{5cm}Quand mon espoir manquait de braise\footnote{\textit{``Chanson du professeur''}, free after G. Brassens}.\newline
\vspace*{0.3cm}
\hspace*{6.0cm}To the memory of Manuel Bronstein
\vspace*{0.3cm}
}

\begin{abstract}
  The first efficient general primality proving method was proposed in the
  year 1980 by Adleman, Pomerance and Rumely and it used Jacobi sums. The
  method was further developed by H. W. Lenstra Jr. and more of his students
  and the resulting primality proving algorithms are often referred to under
  the generic name of \textit{Cyclotomy Primality Proving} (CPP). In the
  present paper we give an overview of the theoretical background and
  implementation specifics of CPP, such as we understand them in the year
  $2007$.
\end{abstract}
\maketitle

\tableofcontents

\section{Introduction}
Let $n$ be an integer about which one wishes a decision, whether it is
prime or not. The decision may be taken by starting from the
definition, thus performing trial division by integers $ \leq
\sqrt{n}$ or is using some related sieve method, when the decision on
a larger set of integers is expected. The method is slow for
relatively small integers, but may be acceptable in certain
contexts. Primality proving becomes a discipline after the realization
that rather than the definition, one may test some property or
consequence of $n$ being prime, and this can often be done by
significantly faster algorithms - basically descending from the
exponential to polynomial asymptotic behavior. This way one easily
eliminates composites which do not verify the particular property of
primes that is tested. The simplest property considered in this
context is certainly \textit{Fermat's small theorem }: $a^{n-1} \equiv
1 \bmod n$ for any $(a, n) = 1$, if $n$ is a prime. Since modular
exponentiation is done in polynomial time in $\log(n)$, such a
\textit{compositeness test} is polynomial.

The disadvantage of the above approach is that there are composites
which verify the same property; such composites are called
\textit{Fermat - pseudoprimes base $a$ } and there is literature
dedicated to these and related pseudoprimes. Stronger statements are
obtained when one has sufficient information about the factorization
of $n-1$. For instance, if there is a prime $q | (n-1)$ and $q >
\sqrt{n}$, while $\left(a^{(n-1)/q}-1,n\right) = 1$ and $a^{n-1}
\equiv 1 \bmod n$, then one easily proves that $n$ is prime.  Indeed,
if $p | n$ is a nontrivial prime factor with $p \leq \sqrt{n}$ -- such
a prime always exists, if $n$ is composite -- then one considers
$\wh{a} = a^{(n-1)/q} \bmod p \in \F_p$. By hypothesis, $\wh{a} \neq
1$ and $\wh{a}^q = 1$; but then $\wh{a} \in \F_p^{\times}$ is an
element of order $q$ and since $\left| \F_p^{\times} \right| = p-1$,
one should have $q | (p-1) < \sqrt{n}$, which contradicts the choice
of $q$. The idea can be refined: $q$ may be replaced by an integer $F
| (n-1), \ F > \sqrt{n}$ which has a known factorization. Based on
this factorization and an easy variation of the above argument, one
obtains a more general primality test. Note that in these cases a
\textit{proof} of primality (or compositeness) comes along with the
result of the algorithm. Tests of this kind can be designed also for
small extensions $\F_{p^k} \supset \F_p$, with astute translations of
the arithmetic in these extensions, in the case when $\F_p$ is
replaced by $\ZM{n}$ and extensions of this ring are used. A general
limitation remains the necessity to know some large factored divisors
$s | (n^k-1)$. Tests of this kind are denoted in general by the name
of \textit{Lucas - Lehmer } tests.

The idea of Adleman et. al. in \cite{APR} was to bypass the above
mentioned restriction, by choosing $k$ so large, that an integer $s >
\sqrt{n}$ and which splits \textit{completely } in small -- albeit, not
polynomial -- prime factors is granted to exist by analytic number
theory. The algebraic part consists in a modification of the Lucas -
Lehmer setting, which allows more efficient testing. In the original 
version of \cite{APR}, the connection to classical test was hard to
recognize. This connection was brought to light by H. W. Lenstra Jr. in
his presentation of the result of Adleman, Pomerance and Rumely at the
Bourbaki Seminar \cite{Le1}.

Let us consider again the Lucas - Lehmer test described above, where
$q | (n-1)$ is a prime with $q > \sqrt{n}$. One can assert that this
test constructs a {\em primitive \nth{q} root } of unity modulo $n$,
in the sense that $\Phi_q(\alpha) = 0 \bmod n$ with $\alpha =
a^{(n-1)/q} \ \rem \ n$ and $\Phi_q(X)$ the \nth{q} cyclotomic
polynomial. It is an important remark, that once $\alpha$ was
calculated, it suffices to verify $\Phi_q(\alpha) \equiv 0 \bmod n$,
and this verification is shorter than the original computation. If $q$
is a proved prime, the verification will yield a proof of primality
for $n$, which can be quickly verified. This is the core idea for
\textit{prime certification}: gathering some information during the
process of an initial primality proof, which can be used for a quicker
a posteriori verification of the proof. Pratt developed this idea in
the context of Lucas - Lehmer tests \cite{Pr}.

When replacing, $q$ by some large factored integer $s$ and searching
for \nth{s} roots of unity $\alpha$ in some extension $\rg{A} \supset
\ZM{n}$, that such roots are zeroes of polynomials over $\ZM{n}$ and
this fact yields a common frame for understanding the APR - test and
generalized Lucas - Lehmer tests. We present here a slight
modification of Lenstra's Theorem 8 in \cite{Le1}, which is seminal to
the approach we take in this paper:
\begin{theorem}
\label{BB}
Let $n > 2$ be an integer and $\rg{A} \supset \ZM{n}$ a commutative
ring extension, $s > 1$, $t = \ord_{s}(n)$ and $\alpha \in
\rg{A}^{\times}$. If the following properties hold
\begin{itemize}
\item[(i)] $\Phi_s(\alpha) = 0$,
\item[(ii)] $\Psi(X) = \prod_{i=1}^t \ \left(X - \alpha^{n^i}\right) \in
\ZM{n}(X)$,
\end{itemize} 
then either $n$ is prime or any divisor $r | n$ verifies:
\begin{eqnarray}
\label{findiv}
  r \in \{ n^i \ \rem \ s : i = 1, 2, \ldots, t = \ord_s(n) \} .
\end{eqnarray}
\end{theorem}
\begin{proof}
Suppose that $n$ is not prime and $r | n$ is a prime divisor. Then
there is a maximal ideal $\eu{R} \supset r \cdot \rg{A}$ which
contains $r$ and $\K = \rg{A}/\eu{R}$ is a finite field while
$\wh{\alpha} = (\alpha \bmod \eu{R}) \in \K$ verifies
$\Phi_s(\wh{\alpha}) = 0$. If $u = \ord_s(r)$, then $\wh{\alpha}^{r^u}
= \wh{\alpha}$ and by galois theory in finite fields, the minimal
polynomial of $\wh{\alpha}$ is $$f(X) = \prod_{i=1}^u \ \left(X -
\wh{\alpha}^{r^i}\right) \in \F_r[X].$$ On the other hand, the
polynomial
\[ \wh{\Psi}(X) = \Psi(X) \bmod \eu{R} = \prod_{i=1}^t \ \left(X -
\wh{\alpha}^{n^i} \right) \in \F_r[X] \] has $\wh{\alpha}$ as zero. By the
minimality of $f(X)$, it follows that $f(X) | \wh{\Psi}(X)$ and since $F_r[X]$
has unique factorization, $\wh{\alpha}^r$ must be a common zero of $f(X)$ and
$\wh{\Psi}(X)$. In particular, there is an exponent $j$ such that
$\wh{\alpha}^r = \wh{\alpha}^{n^j}$ and thus $\wh{\alpha}^{n^j - r} = 1$. But
by (i), $\wh{\alpha}$ is a primitive \nth{s} root of unity, and thus we must
have $n^j - r \equiv 0 \bmod s$ or $r \in < n \bmod s >$. This holds for all
the prime divisors of $n$ and the more general statement \rf{findiv} follows
by multiplicativity.
\end{proof}
This Theorem allows a fundamental generalization of the Lucas - Lehmer tests:
let $n$ be an integer and suppose that an \nth{s} root of unity in the sense
of (i) is found in some ring $\rg{A} \supset \ZM{n}$ and furthermore (ii)
holds. If $s > \sqrt{n}$, then, pending upon a test of the fact that all the
residues
\[ r_i = n^i \ \rem \ s, \quad i = 1, 2, \ldots, t \]
are coprime to $n$, one has a primality proof for $n$. Indeed, if $n$
were composite, then at least one of its prime factors $p \leq
\sqrt{n} < s$. But then, the Theorem implies that $p \in \{ r_i : i =
1, 2, \ldots, t \}$, which is verified to be false. One should note
that prior to Lenstra's work, Lucas - Lehmer tests in ring extensions
of degree $k$ were lacking a transparent criterion for the choice of
the size of the completely factored part $s | (n^k -1)$ required; in
particular, the required factored part was often larger than
$\sqrt{n}$ even for small values of $k$; it was also not possible to
combine informations from tests for different values of $k$ \cite{Wi},
\cite{Mi3}. The Theorem \ref{BB} solves both questions elegantly.

As we have shown in \cite{Mi1}, the Theorem \ref{BB} not only
generalizes the notion of Lucas - Lehmer tests and builds a bridge to
combining them with the test of Adleman, Pomerance and Rumely, it also
indicates a way for a new comprehension of that algorithm. It has
become custom to denote the test described in \cite{APR} in all its
updated variants by \textit{Jacobi sum test}, while \textit{ Cyclotomy
Primality Proving } - or CPP - is a word used to cover all variants of
tests related to Theorem \ref{BB}. These may be in Jacobi sum tests,
Generalized Lucas Lehmer, combinations thereof or also deterministic
variants: we shall indeed see below, that the Jacobi sums test has a
probabilistic \textit{Las Vegas } version, which is mostly the version
used in implementations, and a computationally more complicated
deterministic version. The ideas of CPP were improved by Lenstra
et. al. in \cite{Le2}, \cite{Le3}, \cite{CoLe}, \cite{BH}, \cite{Mi1},
\cite{Mi2},; their constructive base can be described as building a
frame, in which a factor $\Psi(X) | \Phi_s(X) \bmod n$ can be
constructed for some large $s$ and such that, if $n$ is prime, the
factor is irreducible. The computations are performed in {\em
Frobenius rings} extending $\ZM{n}$, which become fields over $\ZM{n}$
if $n$ is prime.

The algorithms of CPP are \textit{de facto} fast, competitive primality
proving algorithms, but they have the complexity theoretical intolerable
feature of a provable {\em superpolynomial} run - time
\begin{eqnarray}
\label{spol}
O\left(\log(n)^{\log \log \log (n)}\right),
\end{eqnarray}
which is the expected value for the size of $t$ in \rf{findiv}. An
practical alternative for proving primality on computers is the random
polynomial test using the group of points of an elliptic curve over
finite fields, originally invented by Goldwasser and Kilian
\cite{GoKi}. The test was made practical by a contribution of
A. O. L. Atkin \cite{AtMo} and has been implemented at the same time
by F. Morain, who maintained and improved \cite{Mo1} a program ECPP
\cite{EC} since more than a decade.

The purpose of this paper is to give a compact presentation of the
theoretical background of the CPP algorithms and an overview of the
basic variants. We also present a new method for computing
\textbf{certificates} of a CPP proof. In the description of
algorithms, we follow a ballance between efficiency and clarity. 

In section two we define the galois, Frobenius and cyclotomic
extensions of rings. The last are the algebraic structures in which
the various tests are performed. Based on this, we then describe an
algorithm for taking roots in cyclotomic ring extensions, which is due
to Huang in the field case. Section three gives an overview of Gauss
and Jacobi sums over galois rings. We then show the connection to the
construction of cyclotomic fields by cyclic field extensions and show
that this mechanism is in fact the core idea of the Jacobi sum
test. In section four we give some computational criteria which
connect this test to the existence and construction of cyclotomic
extensions. In section five we introduce the new certification methods
and the probabilistic algorithms of CPP are defined in section 6.
Finally, in section seven we present the deterministic version of CPP,
and show how it could be understood and implemented as a subcase of
the general CPP test and section 8 contains observations on the run
time and the results from analytic number theory on which the analysis
is based.

The ideas of this paper are updated from the thesis \cite{Mi1} and
many can be found already in the joint thesis of Bosma and van der Hulst 
\cite{BH} and the seminal papers of Lenstra. Our perspective of placing Theorem \ref{BB}
at the center of CPP may be considered as the more personal
contribution of this paper. Based on the common structure of
Lucas - Lehmer and Jacobi sum tests, such as reflected by Theorem
\ref{BB}, we deduce by analogy to the Pratt certificates for Lucas -
Lehmer tests over $\ZM{n}$ a certification method for CPP; such a
method was not known or predicted to exist previously. The same frame
yields also a simple understanding of (a generalized form of) the
Berrizbeitia variant \cite{Be} of the celebrated polynomial time
deterministic test of Agrawal, Kayal and Saxena \cite{AKS}; this is
presented in \cite{AM} and, independently, by Bernstein in \cite{Bern}.

Finally, the notion of cyclotomic extension of rings can be extended
to elliptic extensions of rings - closely connected to the
Schoof-Elkies-Atkins algorithm for counting points on elliptic curves
over finite fields. Together with the use of \textit{dual elliptic
primes}, some relatives of twin primes in imaginary quadratic
extensions of $\Q$, this leads to a new and very efficient combination
of CPP and ECPP (elliptic curve primality proving) algorithms, which
is presented in \cite{Mi4}. The present paper is herewith both an
overview of the recent developments in CPP and a foundation for the
description of new results.
\subsection{Some notations}
Throughout this paper we let $n > 1$ be an integer -- which can be
thought of as a prime candidate. We shall be interested in the ring
$\ZM{n}$ and its extensions and introduce for simplicity the notation
$\id{N} = \ZM{n}$. For integers $s > 0$ we let $\Phi_s(X) \in \Z[X]$
be the \nth{s} cyclotomic polynomial. We shall encounter roots of
unity in various rings. For \textit{complex roots of unity}, we shall
write $\xi_s \in \C$ when $\Phi(\xi_s) = 0$; it will be made clear in
the context, when a certain complex \nth{s} root of unity is fixed. If
$G$ is a finite group and $x \in G$ then $< x >$ will denote the
cyclic group generated by $x$; e.g. $< n \bmod s >$ is the
\textit{cycle} of $n \in \ZM{s}$.  We may at times write
$\log_{(k)}(x)$ for the $k$ - fold iterated logarithm of $x$.  Along
with $n$, we shall often use two parameters $s, t$ in $\N$ such that
$t = \ord_s(n)$ or $s$ is squarefree and $t = \lambda(s) = \lcm_{q |
s} (q-1)$, the product being taken over primes $q$. In both cases, we
consider the following sets related to these parameters:
\begin{eqnarray}
\label{sets}
\id{Q} & = & \left\{ \ q | s \ : \ q \ \hbox{ prime } \right \} \quad
\hbox{ and } \nonumber \\
\id{P} & = & \left\{ \ \wp = ( p^k, q ) \in \N^2 \ : \ p^k || (q-1), \ q
\in \id{Q} \hbox{ and $p$ is prime } \right\}.
\end{eqnarray}
For $\wp = (p^k, q) \in \id{P}$, we may use notations like $p =
p(\wp), k=k(\wp), etc$, with the obvious signification.

\section{Galois extensions of rings and cyclotomy}
Let $\rg{A}$ be a finite commutative ring and $\alpha \in \Omega
\supset \rg{A}$ an element which is annihilated by some polynomial
from $\rg{A}[X]$. Suppose that the powers of $\alpha$ generate a free
module $\rg{R} = \rg{A}[\alpha]$; such modules shall be denoted by
\textit{simple extensions} of $\rg{A}$. Alternately, quotient rings of
the type $\rg{R} = \rg{A}[X]/( f(X) )$, where $f(X) \in \rg{A}[ X ]$
shall also be called simple extensions. It can be verified that
the two types of extensions are equivalent.

There is an ideal $I \subset \Z$ with $I \rg{A} = 0$; the positive
generator $n$ of the annihilator $I$ is the \textit{characteristic} of
the ring $\rg{A}$. We are interested in galois properties of
extensions of finite rings. These have been considered systematically
for primality by Lenstra in \cite{Le1}, \cite{Le2}. The approach we
take here is slightly different and closer to actual computational
aspects; the central concept of \textit{cyclotomic extensions} of
rings end up to be identical to the one of Lenstra.

\begin{definition}
\label{galext}
Let $\rg{A}$ be a finite ring of characteristic $n$ and:
\begin{itemize}
\item[1.] Suppose that there is a galois extension of number fields
  $\KL = \K[ X ]/(f(X))$ with $f(X) \in \id{O}(\K)[X]$, and an ideal
  $\eu{n} \subset \id{O}(\K)$ such that $\rg{A} = \id{O}(\K)/\eu{n}$.
\item[2.] Let $\xi = X + (f(X)) \in \KL, \ \wh{f}(X) = f(X) \bmod
\eu{n} \in \rg{A}[X]$ and
\begin{eqnarray*}
 \rg{R} & = & \id{O}(\KL) / \left( \eu{n} \cdot \id{O}(\LK) \right) = 
\rg{A}[X]/( \wh{f}(X) ) = \rg{A}[\rho], \quad \hbox{ with } \\
 \rho & = & \xi \bmod (\eu{n} \id{O}(\LK) ) = X + (f(X)). 
\end{eqnarray*}
\item[3.] Let $G = \Gal(\KL/\K)$ and for $\sigma \in G$, define
 $\wh{\sigma} : \rho \mapsto \left(\sigma(\xi) \bmod (\eu{n}
 \id{O}(\LK) ) \right)$ and $\wh{G} = \{ \wh{\sigma} : \sigma \in G
 \}$,
\item[4.] Suppose that the degree $d = \deg{f}$, the discriminant
$\disc(f)$ and the characteristic are coprime:
\begin{eqnarray}
\label{units}
\left( \ n, \disc(f) \cdot \deg(f) \ \right) = 1
\end{eqnarray}
\end{itemize} 
If these conditions are fulfilled, then the ring extension $\rg{R}$ is
called \textbf{a galois extension of $\rg{A}$ with group
$\wh{G}$}. Conversely, an extension $\rg{R}/\rg{A}$ is galois, if
there is a galois extension of number fields $\LK/\K$ from which
$\rg{R}$ arises according to 1.-3.
\end{definition}
\begin{remark}
The definition of the galois extension depends in general on the choice of
$(\K, \KL, \eu{n} )$ -- we may in fact identify $\rg{R}$ to this
triple, and unicity of the \textit{lift} to characteristic zero is not a concern. 
In fact, considering the case when $n$ is a prime and $\rg{R} = \F_{n^a}$ is a finite
field, it is obvious that the algebra $\rg{R}$ has multiple lifts. We shall in fact use this
observation and define also when $n$ is not known to be prime, some algebras $\rg{R}$ 
in a simple way, and then construct by operations in $\rg{R}$ additional polynomials 
that split in $\rg{R}$, leading thus to additional lifts to characteristic zero.

The condition \rf{units} is quite artificial, but harmless in
the context of primality testing, where one can think of $\rg{A}$ as
$\id{N}$ or a simple extension thereof: if 4. fails, one has a non
trivial factor of $n$.

\end{remark} 
The main property of a galois extension is of course the fact that the
base ring is fixed by the galois group:
\begin{fact}
\label{galinv}
Let $\rg{R} \supset \rg{A}$ be a galois extension of the finite ring
$\rg{A}$, let $\KL = \K[ X ]/( f(X) )$ be the associated extension of
number fields and $G =\Gal(\KL/\K)$ the galois group. Let $ \rho = X +
\wh{f}(X) \in \rg{R}$ and suppose that $\alpha \in \rg{R}$ is $\hat{G}$ -
invariant. Then $\alpha \in \rg{A}$.
\end{fact}
\begin{proof}
Since $\rg{R} = \pi(\id{O}(\KL))$ -- where $\pi$ is the reduction modulo
$\eu{n} \cdot \id{O}(\KL)$ map -- is a free $\rg{A}$ - module, we can write
$\alpha = \sum_{\sigma \in G} a_{\sigma} \cdot \wh{\sigma} (\rho) \in
\rg{R}$, with $a_{\sigma} \in \rg{A}$. If $\alpha$ is $\hat{G}$ - invariant
and $d = | \hat{G} |  \in \rg{A}^{\ast}$, we have
\[ d \cdot \alpha = \sum_{\sigma, \tau \in G} \wh{\tau} (\alpha) = 
\sum_{\sigma, \tau \in G} a_{\sigma} \cdot \wh{\tau \circ \sigma }( \rho ) =
= A \cdot \Theta, \] with $A = \sum_{\sigma \in G} a_{\sigma} \in
\rg{A}$ and $\Theta = \sum_{\sigma \in G} \wh{\sigma}(\rho) = \pi
\left( \sum_{\sigma \in G} \ \sigma(\xi) \right) = \pi
\left(\Tr_{\KL/\K} \xi \right) \in \rg{A}$. It follows that $\alpha =
(A \cdot \Omega) \cdot d^{-1} \in \rg{A}$, which completes the proof.
\end{proof}
Here are some examples of galois extensions:
\begin{examples}
\label{galex}
\mbox{}
\begin{itemize}
\item[(a)] Let $\rg{A} = \id{N}$, and $s > 0$ such that $(n, s \cdot
\varphi(s)) = 1$. If $f(X) = \Phi_s(X), \K = \Q$ and $\KL =
\Q(\xi_s)$, then $\rg{R} = \Z[\xi_s]/( n \Z[\xi_s] )$ is a galois
extension with group $\hat{G} \sim \ZMs{s}$.

\item[(b)] If $n$ is a prime and $s, \K, \KL, \rg{R}$ are like in the
previous example, let $H \subset \ZMs{s}$ be the decomposition group
of $n$ in $\KL$ with $t = | H |$ and suppose that $(n, s \cdot t ) =
1$. Let $\K_1 = \KL^H$ and $\eu{n} \in \id{O}(\K_1)$ an ideal above
$n$. Then $\id{O}(\K_1)/\eu{n} = \F_n$ and there is a polynomial $f(X)
= \prod_{\tau \in H} \ \left( X - \tau(\xi_s) \right)$ such that
$\rg{r} = \rg{A}[ X ]/(\wh{f}(X)) \subset \rg{R}$ is a galois
extension. Of course, $\rg{r} = \F_{n^t}$ is even a field. 

The construction holds for any subgroup $K$ with $H \subset K \subset
G = \Gal(\KL/\K)$, yielding a filtration of galois extensions. If $K =
G/H$, then for any prime ideal with $n \in \eu{N} \subset
\Z[\zeta_s]$, the group $K$ acts transitively on $\eu{N}$ and $ ( n )
= \prod_{\nu \in K} \nu\left(\eu{N}\right)$.  Furthermore,
$\Z[\zeta_s] / \left(\nu\left(\eu{N}\right)\right) \cong \F_{p^t}$. There
is a canonical decomposition:
\begin{eqnarray}
           \label{candec}
           \rg{R}\  =  \ \Z[\zeta_s]/(n \Z[\zeta_s]) \ \cong \
	      \prod_{\nu \in H} \Z[\zeta_s] /
	      \left(\nu\left(\eu{N}\right)\right) = \prod_{\nu \in H}
	      \ \F_{p^t}.
\end{eqnarray}
If $\rho = \zeta_s \bmod \eu{N} \in \F_{p^t}$ is fixed, then the image
	   of $ \ \zeta_s$ in the Chinese Remainder decomposition of
	   $\rg{R}$ above is $\zeta_s \bmod (n) = \left(\rho,
	   \rho^{k}, \ldots, \rho^{k^{f-1}} \right)$, where $< k \bmod
	   n > = K$ and $\rho^k \equiv \zeta \bmod
	   \sigma^{-1}_{k}(\eu{N})$, with obvious meaning of
	   $\sigma_x(\zeta_s) = \zeta_s^x.\\$

\item[(c)] Let $\id{E}(a,b) : Y^2 = X^3 + aX + b$ be the equation of
an elliptic curve over some field $\K$ and suppose that there is an
ideal $\eu{n} \in \id{O}(\K)$ such that $\id{O}(\K)/\eu{n} = \id{N}$.
Let $\ell$ be a prime and $\psi_{\ell}(X) \in \id{O}(\K)[X]$ be the
\nth{\ell} division polynomial of $\id{E}(a,b)$, its reduction being
$\wh{\psi}_{\ell}(X) \in \id{N}[ X ]$. Then $\rg{R} = \id{N}[X]/\left(
\wh{\psi}_{\ell}(X) \right)$ is a galois extension.
\end{itemize}
\end{examples}

The definition of galois extensions is quite general and is not
specifically bound to the expectation that $n$ might be a prime. We
specialize below galois extensions to Frobenius extensions, which are
related to finite fields.
\begin{definition}
\label{defgalext}
Let $\rg{A}$ be a finite commutative ring of characteristic $n$ and
$\Psi(X) \in \rg{A}[X]$ a monic polynomial. We say that the simple ring
extension $\rg{R} = \rg{A}[X]/(\Psi(X))$ is:
\begin{itemize}
\item[F.] a \textbf{Frobenius} extension, if $\Psi(X^n) \equiv 0 \bmod
\Psi(X)$ and
\begin{itemize}
\item[(F1.)] There is a set $S = \{x_1, x_2, \ldots, x_m \} \subset
\rg{R}$ which generates $\rg{R}$ as a free $\rg{A}$ - module and such that
$\Psi(x_i) = 0$.
\item[(F2.)] There is a group $G \subset \Aut_{\rg{A}} \ (\rg{R} )$
which fixes $S$ and such $g = | G | \in \rg{A}^{\ast}$.
\item[(F3.)] The \textit{traces} of $S$ are $\Tr(x_i) = \sum_{\sigma \in
G} \ \sigma (x_i) \in \rg{A}$.
\end{itemize}
\item[SF.] a \textbf{simple Frobenius} extension, if it is Frobenius
and there is a $t > 0$ such that
\[  \Psi(X) = \prod_{i=1}^t \ \left(X - \zeta^{n^i} \right), \quad 
\hbox{ where } \quad \zeta = X + (\Psi(X)) \in \rg{R}. \]
\end{itemize}
\end{definition}
\begin{remark}
\label{crux}
The example (c) is a galois extension which is in general not
Frobenius. The other two examples are Frobenius at the same time and
the first extension in (b) is simple Frobenius. The property
F3. implies that $\rg{A}$ is exactly the ring fixed by $G$, the proof
being similar to the one of Fact \ref{galinv}.

The situation in (b) is crucial for CPP. In fact, in the algorithms we
shall investigate integers $n$ that lead to the decomposition
\rf{candec} of the ring $\Z[\zeta_s]/(n \Z[\zeta_s])$ and show that
for such integers the Theorem \ref{BB} can be applied.
\end{remark}
Next we clarify the notion of \textit{primitive root of unity}, which
has some ambiguity when considering roots of unity over rings. The
question is illustrated by the simple example: is $a = 4 \bmod 15$ a
primitive second root of unity modulo $15$ ? One verifies that $a \in
\ZMs{15}$ and $a \neq 1$ while $a^2 = 1$. However, $a-1 \not \in
\ZMs{15}$. We shall avoid such occurrences and define
\begin{definition}
\label{pru}
Let $\rg{A}$ be a commutative ring with $1$ and $s > 1$ an integer,
while $\Phi_s(X) \in \Z[X]$ is the \nth{s} cyclotomic polynomial,
$\wh{\Phi}_s(X) \in \rg{A}[X]$ its image over $\rg{A}$. We say that
$\zeta \in \rg{A}$ is a \textit{primitive \nth{s} root of unity} iff
$\wh{\Phi}_s(\zeta) = 0$.
\end{definition}

We leave it as an exercise to the reader to verify that, if $\rg{A}$
  is finite, then $\zeta \in \rg{A}$ is an \nth{s} primitive root of
  unity if and only if for all maximal ideals $\eu{A} \subset \rg{A}$,
  the root of unity $\zeta_{\eu{A}} = \left( \zeta \bmod \eu{A}
  \right) \in \K = \rg{A}/\eu{A}$ is a primitive \nth{s} root of unity
  in the field $\K$. In particular, $\zeta - 1 \in
  \rg{A}^{\times}$. Note that $\zeta_s \bmod ( n \Z[\zeta_s] )$ in
  example (b) is a primitive root of unity.

The next step towards the goal of Remark \ref{crux} consists in
  defining the cyclotomic extensions of rings, which are simple
  Frobenius extensions generated by primitive roots of unity.

\begin{definition}
\label{cycext}
Let $n, s$ and $\id{N}$ be as above and $\Omega \supset \id{N}$ some
ring with $\zeta \in \Omega$, a primitive \nth{s} root of unity. We
say that
\[ \rg{R} = \id{N}[ \zeta ] \] 
is an \nth{s} \textit{cyclotomic extension of the ring } $\id{N}$, if
the extension $\rg{R} / \id{N}$ is simple Frobenius. In particular,
$\rg{R}/\id{N}$ has the galois group $G = < \sigma > $ generated by
the automorphism with $\sigma(\zeta) = \zeta^n$ and $| G | = t =
\ord_s(n)$.

We say that $s$ is the {\em order} and $t$ is the {\em degree}
of the extension ${\bf R}$. Sometimes we shall denote the extension
also by the triple $(\rg{R}, \zeta, \sigma)$.
\end{definition}

Like for finite fields, a galois extension $\rg{R} \supset \id{N}$ can
be an \nth{m} cyclotomic extension of $\id{N}$ for various values of
$m$. We shall in fact often start with galois extensions $\rg{R}$ of
degree $d$ over $\id{N}$ and then seek \nth{m} primitive roots of
unity in $\rg{R}$, for various values $m | (n^d - 1)$ and then prove
that these roots together with the galois group generate an \nth{m}
cyclotomic extension. The procedure will be illustrated below, in the
results on the Lucas -- Lehmer test.

It is also natural to consider \textit{subextensions} of cyclotomic
extensions, i.e. rings of the kind
\[ \rg{T} = \id{N}[ \eta ], \quad \hbox{ with } \quad \eta = \sum_{i=1}^{t/u} 
\sigma^{u i} (\zeta) \in \rg{R} ,\] where $u | t$. Such subextensions
are galois (even abelian). They have been considered recently by
Lenstra and Pomerance in their version of the AKS algorithm \cite{LP}; the term
of \textit{pseudo - fields } was coined in that context.
\begin{remark}
\label{subext}
Let $(\rg{R}, \zeta, \sigma)$ be some \nth{s} cyclotomic extension of
$\id{N}$, with $\rg{R} = \id{N}[\zeta]$ and $t = [ \rg{R} : \id{N} ]$. Suppose
that there is an integer $u > 1$, and $\beta \in \rg{R}$ with $\Phi_u(\beta) =
0$ and such that $\rg{S} = \id{N}[\beta]$ is an \nth{u} cyclotomic extension
with automorphism group induced by the restriction of $\sigma$ to $\rg{S}$. We
shall say in such a case, by abuse of language, that $(\rg{R}, \beta, \sigma)$
is a \nth{u} cyclotomic extension. 
\end{remark}

Cyclotomic extensions do not exist for any pair $(n, s)$ and their existence
is a (sporadic) property of the number $n$ with respect to $s$; this fact is
used for primality testing. The following theorem groups a list of equivalent
properties of cyclotomic extensions, relating them to Theorem \ref{BB} and
providing a useful base for algorithmic applications.

\begin{theorem}  
\label{th1}
Let $s, n > 1$ be coprime integers, $t = \ord_s(n)$ be also coprime to
$n$, and fix $\xi_s \in \C$. Let $\rg{A}$ be the ring of integers in
$\KL_n = \Q[\xi_s]^{<n \smod s>}$ and consider the polynomial
$\Psi_0(x) = \prod_{i=1}^{t-1} \ (x - \xi_s^{n^i}) \in \rg{A}[x]$. The
following statements are equivalent:
\begin{itemize}
\item[(I)] An \nth{s}  cyclotomic extension of  $\id{N}$  exists.

\item[(II)] $r | n \ \Longrightarrow r  \ \in \ <n  \mod s >$.  

\item[(III)] There is a surjective ring homomorphism $\tau_0 : \rg{A}
             \rightarrow \id{N}$.
\item[(IV)] There is a polynomial $\Psi(x) = \tau_0 \left(\Psi_0
\right) \in \id{N}[x]$ of degree $t$ with:
\begin{itemize}
\item[(i)] $\Psi(x) \, | \, \Phi_s(x)$ 

\item[(ii)] if $\zeta = x + (\Psi(x)) \in \id{N}[x] / (\Psi(x))$, 
            then $\Psi(\zeta^{n^i}) = 0$, for $i=1,2,\ldots,t. $
\end{itemize}
\end{itemize}
\end{theorem}
\begin{proof}
Suppose that $(I)$ holds and let $\Psi(X) = \prod_{i = 1}^t \left(X -
\sigma^i(\zeta) \right) \in \id{N}[X]$. The argument used in the proof
of Theorem \ref{BB} shows that $(I) \Longrightarrow (II)$.

Assume that $(II)$ is verified, $r | n$ be some prime factor and let
$\rho \in \overline{\F}_r$ be a primitive \nth{s} root of unity. If
$\eu{R} \subset \KL_r = \Q(\xi_s)^{< r \bmod s >}$ is some prime ideal
above $r$, then $\id{O}(\KL_r) \bmod \eu{R} = \F_r$ as follows from
the example (c) and relation \rf{candec}. But since $r \in < n \bmod s
>$ it follows that $\rg{A} \subset \id{O}(\KL_r)$ and there is a
fortiori a surjective map $\tau_r : \rg{A} \rightarrow \F_r$. By
Hensel's Lemma and the Chinese Remainder Theorem, this map can then be
extended to a map $\tau_0 : \rg{A} \rightarrow \id{N}$, so $(II)
\Rightarrow (III)$.

Assume $(III)$ holds and let $r | n$ be a prime. Then $\tau_0$ extends by
composition with the reduction modulo $r$ to a map $\tau_r : \rg{A}
\rightarrow \F_r$. In particular $\Psi_r(X) = \tau_r\left(\Psi_0(X)
\right) \in \F_r[ X ]$ is a polynomial such that $\Psi_r(X) |
\Psi_s(X)$ and $\Psi_r(\zeta^n) = 0$ if $\Psi_r(\zeta) = 0$. Using
again Hensel's Lemma and the Chinese Remainder Theorem, a polynomial
$\Psi(X) \in \id{N}(X)$ with the same properties can be constructed
and thus $(III) \Rightarrow (IV)$
 
Finally, if $\Psi(x) \in \id{N}[x]$ has property $(IV)$, let $\rg{R} =
\id{N}[x] / (\Psi(x))$ and $\zeta = x + (\Psi(x))$; it follows from
(i) that $\zeta$ is a primitive \nth{s} root of unity. We have to show
that $\sigma: \zeta \mapsto \zeta^n$ is an automorphism of
$\rg{R}$. By construction, $\sigma$ permutes the zeroes of $\Psi$, so
$G = < \sigma > $ acts transitively on $S = \{ \zeta, \zeta^n, \ldots,
\zeta^{n^{t-1}} \}$. This shows that $\rg{R}$ is cyclic Frobenius, and
since $\Psi(X) | \Phi_s(X)$, it is an \nth{s} cyclotomic extension of
$\id{N}$, so $(IV) \Rightarrow (I)$.
\end{proof}

\begin{remark}
\label{corerem}
It follows from $(III)$, that the extension $\rg{R}/\id{N}$ is galois
in the sense of the Definition \ref{galext} and this confirms the fact
that its subextensions are galois too.

The relation (ii) is an elementary verifiable condition for the
existence of cyclotomic extensions. If $\lambda(s) = \lcm_{(q = p^e)
|| n} (\varphi(q) )$ is the Carmichael function, while $\varphi$ is
Euler's totient function, then $\ZMs{s}$ contains $\varpi(s) =
\frac{\varphi(s)}{\lambda(s)}$ disjoint cyclic subgroups. The larger
$\varpi(s)$, the more improbable it becomes to find integers $n$ for
which (ii) is verified. This is the core idea of the CPP tests.
\end{remark}

The following simple fact has some important implications about the
size of cyclotomic extensions.

\begin{fact}
Let $p$ be a prime, $n \in \N_{>1}, (n, p) = 1$ and $v_p(x), \ x \in
\Z$ denote the $p$ - adic valuation. If $p$ is odd, $t = \ord_p(n)$
and $v = v_p(n^t - 1)$, then the order
\[ \ord_{p^m}(n) = t \cdot p^u \quad \hbox{ with } \quad u = \max(0, m - v) .\]
For $p = 2$ we distinguish the following cases:
\begin{itemize}
\item[1.] If $n \equiv 1 \bmod 4$ and $v = v_2(n-1)$, then
\[ \ord_{2^m}(n) = 2^u \quad \hbox{ with } \quad u = \max(0, m-v)  ,\]
\item[2.] If $n \equiv 3 \bmod 4$ and $v = v_2(n+1)$, then 
\[ \ord_{2^m}(n^2) = 2^{u+1} \quad \hbox{ with } \quad u = \max(0, m - v). \]
\end{itemize}
\end{fact}
\begin{proof}
The proof is left as an exercise to the reader, see also \cite{Se}, Chapter II, \S 3.
\end{proof}
The remarkable phenomenon above consists in the fact that the order
$\ord_{p^m}(n)$ starts from an initial value $t = \ord_p(n)$ which is
constant for $m \leq v_p(n^t-1)$ and then increases by factors $p$,
when $m$ grows. The only exception is the case $p = 2$ and $n \equiv 3
\bmod 4$, when one has to consider $m \leq v_2(n^2 - 1)$ as starting
value. This leads to the following
\begin{definition}
\label{satext}
Let $n$ be an integer and $p$ a prime with $(p, n) = 1$. 
We define the \textit{saturation index} of $n$ with respect to $p$ by:
\begin{eqnarray}
\label{satind}
        k_n(p) & = & \left\{ \begin{array}{ll} v_2(n^2-1) & \hbox{if
        $ \ p = 2 $ and $ \ n \equiv 3 \mod 4, $ } \\ & \\ v_p(n^t-1) &
        \hbox{with} \quad t = \ord_p(n) \quad \hbox{ otherwise. }
        \end{array} \right.
\end{eqnarray}
If $(\rg{R}, \sigma, \zeta)$ is a \nth{p^k} cyclotomic extension of
$\id{N}$ and $k \geq k_n(p)$, then the extension is \textit{saturated
}. In general, if $(s, n) = 1$ and $(\rg{R}, \sigma, \zeta)$ is an
\nth{s} cyclotomic extension of $\id{N}$, we say that the extension is
saturated, if $p | s \Rightarrow p^{k_n(p)} | s$.

For odd $p$ or $p = 2$ and $n \equiv 1 \bmod 4$, we shall denote by
\textit{saturated \nth{p}} extension of $\id{N}$ a galois extension
with $[ \rg{R} : \id{N} ] = d, \ (d, p) = 1$ or $d = 2$, if $p = 2$
and $n \equiv 3 \bmod 4$, and which is a \nth{q} cyclotomic extension
of $\id{N}$, with $q = p^{k_n(p)}$.
\end{definition}

Note that the term {\em saturated \nth{p^h} extension}, implicitly
asserts the fact that $h \geq k_s(p)$; the definition of a {\sl saturated
\nth{p}} cyclotomic extension is an exception, since it denotes an
extension which not only contains a \nth{p} root but also a
\nth{p^{k_n(p)}} primitive root of unity. It can happen that a \nth{p}
cyclotomic of $\id{N}$ exists, but not a \textit{saturated} one, as
illustrated by:

\begin{example}
\label{nonsat}
Let $n = 91 = 7 \cdot 13$. Then (II) implies that a third cyclotomic
extension of $n$ exists, since $r \in \ <n \mod 3>$ for all $r|n$.
However, according to $(II)$ of Theorem \ref{th1}, this extension is
not saturated, since $n = 1 \mod 9$ yet $r = 7 | n$ and $r \not \in
\left( \ < n \bmod 9 > \ = \ < 1 > \ \right)$.
\end{example}

The saturated extensions are characterized by the following property:

\begin{theorem}
\label{thsat}
If $(\rg{R},\sigma,\zeta)$ is a saturated \nth{p^{k_s}} extension and
$h > k_s$ then a \nth{p^h} extension of $\id{N}$ exists.
\end{theorem}
\begin{proof}
Consider $h>k_s$, $\rg{R}(h) = \rg{R}[x]/(x^{p^{(h-k_s)}} - \zeta)$
and let $\zeta(h)$ be the image of $x$ in~$\rg{R}(h)$. It is easy to
establish by comparing ranks, that $(\rg{R}(h),\sigma(h),\zeta(h))$ is
a \nth{p^h} cyclotomic extension -- where $\sigma(h)$ is the
extension of $\sigma$ to $\rg{R}(h)$.
\end{proof}
Theorem \ref{thsat} motivates the denomination of ``saturated'': the
existence of a saturated \nth{p} extension implies existence of
cyclotomic extensions of degree equal to any power of $p$. The Example
\ref{nonsat} shows that the existence of a saturated extension is also
necessary for this. We shall use for commodity, the term of {\sl
complete extension} for the union of all saturated extensions of orders
$p^h$:
\begin{definition}
\label{compsat}
Suppose that a saturated \nth{p^{k}} extension $(\rg{R},\sigma,\zeta)$
of $\id{N}$ exists and let:
\begin{eqnarray}
\label{satexp}
         \bigl(\rg{R}_{\infty} (p),\sigma_{\infty} (p), \zeta_{\infty}
                (p) \bigr) & = & \bigcup_{h = k_s}^{\infty} \
                (\rg{R}(h),\sigma(h),\zeta(h)).
\end{eqnarray}
$\rg{R}_{\infty}$ is called {\it complete} \nth{p} extension and 
its existence is granted by the premises and the preceding theorem.
\end{definition}

The proving of existence of cyclotomic extensions focuses herewith on
proving existence of saturated extensions. The existence of saturated
extensions has also implications for the properties of primes $r$
dividing $n$ \cite{CoLe}:

\begin{lemma}[Cohen and Lenstra, \cite{CoLe}]
\label{lp}
Suppose that $p$ is a prime with $(p,n) = 1$, for which a saturated
\nth{p} cyclotomic extensions of $\id{N}$ exists. Then for any $r | n$
there is a $p$-adic integer $l_p(r)$ and, for $p>2$, a number $u_p(r)
\in \ZM{(p-1)}$, such that:
\begin{eqnarray}
\label{lpr}
r & = & n^{u_p(r)} \mod p \quad \mbox{and} \nonumber \\ & & \nonumber
  \\[-0.3cm] r^{p-1} & = & (n^{p-1})^{l_p(r)} \ \in \ \{1+p \cdot \Z_p
  \} \ \ \ \mbox{if} \quad p > 2, \\ & & \nonumber \\ [-0.3cm] r & = &
  n^{l_p(r)} \ \in \ \{1+2 \cdot \Z_2\} \quad \mbox{if} \quad p =
  2.\nonumber
\end{eqnarray}
\end{lemma}

\begin{proof}
Using Theorem \ref{thsat}, the hypothesis implies that $r \in \ <n
\mod p^k>$ for all $k \geq 1$ which implies \rf{lpr}.
\end{proof}

\subsection{Finding Roots in Cyclotomic Extensions}
Consider the following problem: given a finite field $\F_q =
\F_p[\alpha]$ with $q = p^k$ a prime power and $r$ a prime with
$v_r(q-1) = a$, and given $x \in \F_q$ with $x^{(q-1)/r} = 1$, find a
solution of the equation $y^r = x$ in $\F_q$. The problem has an
efficient polynomial time solution, if a \nth{r^a} root of unity $\rho
\in \F_q$ is known and the algorithm was described by Huang in
\cite{Hu1, Hu2}.

We shall treat here the generalization of the problem to cyclotomic
extensions of rings. The basic idea is the same and it is well
illustrated by the case $r = 2$ and $q = p \equiv 5 \bmod 8$. In this
case we let $u = x^{(p-1)/4} = \pm 1$, since $u^2 = x^{(p-1)/2} = 1$
by hypothesis. But $e = (p-1)/4$ is odd and thus $f = (e+1)/2$ is an
integer, while $x^{2f} = u \cdot x$. If $\rho^2 = u = \pm 1$ for $\rho
\in \F_p$, then a solution of $y^2 = x$ is given by $y = \rho^{-1}
\cdot x^f$. Thus, knowing a \nth{4} root of unity, one can find
square roots in $\F_p$. The general case is described in the following:

\begin{theorem}
\label{troot}
Let $p$ be a prime with $(p, n) = 1$ and $(\rg{R}, \sigma, \zeta)$ a
saturated $p$-th cyclotomic extension of $\id{N}$; let $\alpha \in
\rg{R}$ and $l \leq k_p(n)$ be such that
\begin{equation}
\label{pow}
\alpha^{N/p^l} = 1
\end{equation}
is satisfied. Then there is a polynomial deterministic algorithm for
finding a root $\beta \in \rg{R}$ of the equation $x^{p^l} = \alpha$.
\end{theorem}
\begin{proof}
Let
\[ t = \sharp < \sigma > = [ \rg{R} : \id{N} ], \quad N = n^t - 1, \quad k
 = v_p(N) , \] and let $u$ be given by $N = u \cdot p^k$, so that $(u,
 p) = 1$ and $k = k_n(p)$. Since \rg{R} is saturated, $\zeta \in
 \rg{R}$ is a \nth{p^k} root of unity. If $\alpha$ is a \nth{p^l}
 power in \rg{R}, then
\begin{equation}
\label{lpow}
	\alpha^u = \zeta^{\nu \cdot p^l}, \quad \hbox{with} \quad \nu
	\in \ZM{p^{k-l}} .
\end{equation}
Note that $\nu \mod p^i$ can be successively computed for $i = 1, 2,
\ldots, k-l$ by comparing $\alpha^{u \cdot p^{k-i-l}}$ to powers of
$\zeta^{p^{k-i}}$. Given $\nu$, one can define a solution to $x^{p^l}
= \alpha$ in the following way. Let $u'$ be such that $u \cdot u' = -1
\mod p^k$ and $e = (1+u \cdot u') /p^l$. Then $\beta = \alpha^e \cdot
\zeta^{-u' \cdot \nu}$ is such that $\beta^{p^l} = \alpha$, which
follows from a straightforward computation.
\end{proof}

\subsection{Finding Roots of Unity and the Lucas -- Lehmer Test}
The algorithm described above assumes that a saturated root of unity
is known in a galois extension of appropriate degree. This can be
found naturally by trial and error. Suppose that one wants to
construct a saturated \nth{p} cyclotomic extension and $t =
\ord_p(n)$. The \textit{bootstrapping} problem that one faces,
consists in finding first a galois extension $\rg{R}/\id{N}$ of degree
$t$; if such an extension is provided, one seeks a \nth{p} power non
residue, like one would do if $\rg{R}$ was a field.

Let us recall some facts and usual notations about cyclotomic fields (see also
\cite{Wa}). The \nth{s} cyclotomic field is $\LK_s = \Q(\zeta_s) = \Q[ X ]/ (
\Phi_s(X) )$, an abelian extension of degree $\varphi(s)$ with ring of
integers $\id{O}(\KL_s) = \Z[\zeta_s]$ and galois group 
\[ G_s = \Gal(\KL_s/\Q) = \left\{ \sigma_a : \zeta_s \mapsto \zeta_s^a ;
\hbox{ where } \ (a, s) = 1 \right\} \ \cong \ \ZMs{s} .\] It is noted
in \cite{Le2} that in fact $\sigma_a = \lchooses{\KL_s/\Q}{a}$ is in
this case the Artin symbol of $a$. We shall adopt the notation of
Washington, introduced above. The Theorem of Kronecker - Weber states
that all abelian extensions of $\Q$ are subfields of cyclotomic
extensions and if $\K/\Q$ is an abelian field, then its
\textit{conductor} is by definition the smallest integer $s$ such that
$\K \subset \KL_s$, with $\KL_s$ the \nth{s} cyclotomic extension.

The next fact shows where to look for galois extensions of $\id{N}$.
\begin{fact}
\label{exts}
Let $n > 2$ be an integer and $\K /\Q$ be an abelian extension of conductor
$s$ such that $(s, n) = 1$ and $\Gal(\K/\Q) \cong < n \bmod s >$, $t =
\ord_s(n) = [ \K : \Q ]$. Then there are $\omega_i \in \id{O}(\K), \ i = 1, 2,
\ldots, t$ such that
\begin{eqnarray}
\label{repr}
 \id{O}(\K) = \Z[\omega_1, \omega_2, \ldots, \omega_t ] \quad \hbox{ and }
 \quad \rg{R} = \id{N}[ \wh{\omega}_1, \wh{\omega}_2, \ldots, \wh{\omega}_t ], 
\end{eqnarray}
where $\wh{\omega}_i = \omega_i \bmod n \id{O}(\K)$.

The ring $\rg{R} = \id{O}(\K)/\left(n \cdot \id{O}(\K) \right)$ is a galois
extension of $\id{N}$ if and only if the ring $\id{O}(\K)$ has a normal $\Z$ -
base. 

If $t = p^k$ is a prime power, this happens in the following cases:
\begin{itemize}
\item[(i)]    $s$ is prime and $p^k \parallel (s-1)$.
\item[(ii)]   $p$ is odd and $s=p^{k+1}$.
\item[(iii)] $p = 2, \ k \geq 2$ and $s = 2^{k+2}$. 
\item[(iv)] $p = 2$ and $k = 1$.
\end{itemize}
\end{fact}
\begin{proof}
The ring $\id{O}(\K)$ is a free $\Z$ - module of rank $t$ and discriminant
which divides $s^N$ for some integer $N > 1$, see e.g. \cite{Ri}. With this,
the assertions become simple verifications based upon the definition of $\K$ and
the one of a galois extension. The assumption that $t = p^k$ can be dropped, by
using the linear independence of \nth{s_i} cyclotomic extensions ($i = 1, 2$)
when $(s_1, s_2) = 1$.
\end{proof}

The following Theorem is useful for constructing roots of unity and
the associated cyclotomic extensions, as well as for generalized
Lucas-Lehmer tests:
\begin{theorem}
\label{LL}
Let $\id{N}$ and $\K$ be a field of conductor $s$ as described in Fact
\ref{exts}; in particular $\rg{R} = \id{O}(\K)/(n \cdot \id{O}(\K))$ is
a galois extension of $\id{N}$ and with $\Gal(\KL/\Q) = \ <n \mod
s>$. Suppose that
\begin{eqnarray}
\label{autom}
\exists \ \alpha \in O(\KL) & \mbox{such that} & \sigma_n(\alpha) =
\alpha^n \mod n \cdot O(\KL).
\end{eqnarray}
For all primes $q|s$, let $(n^t-1)/q = \displaystyle\sum_{i=0}^{t-1} \
c_i(q) \cdot n^i$ and suppose that
\begin{eqnarray}
\label{pows}
\beta(q)  =  \prod_{i=0}^{t-1} \ (\sigma_n^i(\alpha^{c_i(q)})) \mod
n \cdot O(\K) \quad \hbox{ verifies } \quad (\beta(q) - 1) \in
\rg{R}^{\ast}.
\end{eqnarray}
Let $\beta =\prod_{q|s} \beta(q) \mod n \cdot O(\K)$ and $\sigma$ be the
automorphism induced by $\sigma_n$ in $\rg{R}$. Then $(\rg{R},\sigma,\beta)$
is a {\em saturated } \nth{s} cyclotomic extension of $\id{N}$.
\end{theorem}
\begin{proof}
Let $ \Psi(X) \ = \ \prod_{i=0}^{t-1} \ (X - \sigma(\beta))$; then $\Psi(X)
\in \id{N}[X]$ since $\sigma_n$ generates $\Gal(\K/\Q)$. Furthermore $\Psi(X)
| \Phi_s(X)$ by construction and thus $\beta$ is a primitive \nth{s} root of
unity. The statement follows from $(IV)$ of Theorem \ref{th1}. Note that
$\id{N}[\beta(q)]$ are by construction saturated extensions, and thus
$\id{N}[\beta]$ is saturated too.
\end{proof}
\begin{remark}
\label{LLapp}
In practical applications of Theorem \ref{LL}, $n \gg q$. We show that
there is a simple expansion of the shape $(n^t - 1)/q =
\sum_{i=0}^{t-1} \ c_i(q) \cdot n^i$ which makes the computation of
$\beta(q)$ in \rf{pows} particularly efficient. Let $n = a \cdot q +
b$, with $0 \leq b < q$. Then
$$(n^t-1)/q \ = \  ((n^t - b^t) + (b^t-1))/q \ = \  
(b^t-1)/q+a \cdot \sum_{i=0}^{t-1} \ n^{t-i-1} \cdot b^i.$$ This leads to the 
equation:
\begin{eqnarray}
\label{fastpow}
(n^t-1)/q \ = \ \sum_{i=0}^{t-1} \ c_i(q) \cdot n^i, & & \mbox{with} \
c_0(q) = (b^t-1)/q \\ 
               & & \mbox{and} \ c_i(q) = a \cdot b^i, \ \mbox{for} \
                i=1,\ldots,t-1. \nonumber
\end{eqnarray}
Note that $c_i(q)$ are not necessarily $<n$, but $c_i(q)/n$ is at most
a small number. The regular shape of the coefficients is however very
useful for the simultaneous computation of $\alpha^n$ and
$\alpha^{c_i(q)}$.  Suppose that the cost for the application of one
automorphism $\sigma_n$ is $c \cdot \mbox{(multiplication in
$\rg{R}$)}$ -- if no fast polynomial multiplication methods are used,
then $c = 1$. The time needed for the evaluation of $\beta(q)$ using
\rf{fastpow} is bounded by
\[ 2 \cdot (t-1)(\log q + c + 1)+\log n. \]
This method of evaluation is thus 
for $t > (\log q)/\log(n / q^2 \cdot 4^{c+1})$ more efficient 
than when defining $\beta(q) = \alpha^{(n^t-1)/q}$, and performing the direct 
exponentiation.   
\end{remark}

\NI If $\K$ is a field of degree $p^k$ defined by Fact \ref{exts}, an
\nth{s} cyclotomic extension can be constructed by using Theorem
\ref{LL}. This is the Lucas - Lehmer approach to constructing
cyclotomic extensions. It is obvious that, when the degree of
extensions is of importance and the order irrelevant, a minimal $s$
will be chosen.
\begin{remark}
The extensions constructed by the Lucas - Lehmer method are {\sl
saturated}. This approach is used in \cite{AdLe} for constructing
galois fields. We shall also show that this has useful consequences
for combining cyclotomic extensions.
\end{remark}

\section{Gauss and Jacobi sums over Cyclotomic Extensions of Rings}
Gauss Sums are character sums used in various contexts of
mathematics. It will be important for us to note that Gauss sums are
Lagrange resolvents encountered when solving the equation $X^s = 1$
with radicals, over $\Q$. Or, equivalently, when building the \nth{s}
cyclotomic field $\KL_s/\Q$ by a succession of prime power galois
extensions, see e.g. \cite{La}.

Let $n, m > 1$ be integers with $m$ squarefree and let $\lambda(m)$
be the exponent of $\ZMs{m}$, where $\lambda$ is the Carmichael
function. In this section, $u = \lambda(m)$ and $f$ will be some
divisor of $u$ and we assume that $(n, m u) = 1$. Let $\rg{A} \supset
\id{N}$ be a galois extension which contains two primitive roots of
unity $\zeta, \rho$ of respective orders $u, m$. A multiplicative
character $\chi: \ZMs{m} \rightarrow \ < \zeta >$ is a multiplicative
group homomorphism $\ZMs{m} \rightarrow \ < \zeta >$.  We denote by
$\ZMd{m}$ the set of multiplicative characters defined on $\ZMs{m}$;
the set $\ZMd{m}$ builds a multiplicative group and the {\em order} of
$\chi$ is the cardinality of the image $\mbox{Im}(\chi)$. We shall
also denote characters $\chi \in \ZMd{m}$ by characters modulo $m$.

Let $\chi \in \ZMd{m}$ and $d$ be a divisor of $m$. If there is a
character 
\[\chi': \Z \ \Big/ \left({\frac{m}{d}} \cdot \Z \right)
\rightarrow \ <\zeta> \] 
such that
\[ \chi(x) = \chi'\left(x \bmod (m/d) \right) \quad \hbox{for all} \quad x \in
\ZMs{m}, \] 
then $\chi$ is said to be {\sl induced} by $\chi'$. A character $\chi:
 \ZMs{m} \rightarrow \ < \zeta >$ is called {\sl primitive} if it is
 induced by no character different from itself; in this case, $m$ is
 called the {\sl conductor} of $\chi$. Each character $\chi$ is
 induced by a unique primitive character $\chi'$ and the conductor of
 $\chi$ is defined to be equal to the conductor of the primitive
 character it is induced by. In particular, the principal character
 ${\bf 1}: \{1\} \rightarrow \ <1>$ is primitive and has conductor 1.

For $\chi \in \ZMd{m}$, we shall set for ease of notation 
\[ \chi(x) = 0 \quad \hbox{for } \quad (x, m) > 1.\]
 The {\sl Gauss-Sum} of $\chi$ with respect to $x$ is the element of
\rg{A} given by\footnote{We adopt Lang's sign definition for the character sums}:
\begin{eqnarray}
\label{2.18}
	\tau(\chi) \ = -\sum_{x \in \ZM{m}} \ \chi(x) \cdot \rho^x .
\end{eqnarray}
The Gauss-Sum depends upon the choice of an element in $<\rho>$
according to:
\begin{eqnarray}
\label{2.19}
	\tau_a(\chi) \ = -\sum_{x \in \ZM{m}} \ \chi(x) \cdot \rho^{a
	 \cdot x} \ = \ \chi^{-1} (a) \cdot \tau(\chi), \ \ \ \
	 \forall \, a \in \ZMs{m}.
\end{eqnarray}
Let $\nu \in \ZMs{m}, \ H(\nu) = \ZMs{m} / < \nu \mod m>$ and $h$ be a
coset in $H( \nu )$. The \nth{h} {\sl Gauss Period} with respect to
$\nu$ is defined by:
\begin{eqnarray}
\label{2.20}
	\eta_h(\rho, \nu) \ = \ \sum_{\mu \in h} \ \rho^{\mu} , \ \ \ \
	\forall \, h \in H(\nu).           
\end{eqnarray}
Let $H( \nu )^{\widehat{} \ } = \{ \ \chi \in \ZMd{m} \mid \chi( \nu )
= 1 \ \}$.  $H(\nu)^{\widehat{} \ }$ is dual to $H( \nu )$ in the
sense that characters $\chi \in H(\nu)^{\widehat{} \ }$ operate on
cosets $h \in H(\nu)$.  Gauss-Sums and Gauss Periods are connected by:
\begin{eqnarray}
\label{2.21}
 \tau(\chi) = -\sum_{h \in H(\nu)} \ \chi(h) \cdot \eta_h(\rho, \nu)
\end{eqnarray}
and 
\begin{eqnarray}
\label{2.22}
	\mid H( \nu ) \mid \ \cdot \ \eta_h(\rho, \nu) \ = \ 
        -\sum_{\chi \in H(\nu)^{\widehat{} \ }} \ 
		\chi^{-1}(h) \cdot \tau(\chi),\ \ \ \ \ \
	\forall \, h \in H(\nu).           
\end{eqnarray}
Equation \rf{2.21} follows by using $\chi( \nu ) = 1$ and regrouping
the summation order in \rf{2.18}. Identity \rf{2.22} is a consequence
of the following:

\begin{fact}
\label{2.23}
If  $G$  is a subgroup of  $\ZMd{m}$  and  $x \in \ZMs{m}$, then
\begin{eqnarray}
\label{2.24}
	s(x) \ = \ \sum_{\chi \in G} \ \chi(x) & = &
		\left\{\begin{array}{ll} 0, & \hbox{if $\ \ \exists
		\chi \in G$ with $\chi (x) \neq 1$,} \\ & \\ \mid G
		\mid , &\hbox{otherwise.}
\end{array} 
\right. 
\end{eqnarray}
\end{fact}
\begin{proof}
If $\chi(y) \neq 0$, then $s(x) \cdot (1 - \chi(y)) = 0 $ leads to the
claimed result.
\end{proof}

If $\chi, \chi' \in \ZMd{m}$ are primitive, the {\it Jacobi-Sum} $ j(
\chi, \chi')$ is defined by:
\begin{eqnarray}
\label{2.25}
 j( \chi, \chi') \ = -\sum_{x \in \ZM{m}} \
	\chi(x) \cdot \chi' (1-x).
\end{eqnarray}
Gauss and Jacobi-Sums factor with respect to the ideals $\bigl(
p^{v(p)} \cdot \ZM{m} \bigr)$ of $\ZM{m}$, where $p^{v(p)} \parallel
m$ and this will allow us to restrict our attention to characters of
prime conductor. The factorization is given by the following:
\begin{fact}
\label{2.26}
Let $\chi \in \ZMd{m}$ and $m = \prod_{p \mid m} \ p^{v(p)}$.  Then
there are characters $\chi_p \in \ZMd{p^{v(p)}}$ and Gauss-Sums
$\tau_p(\chi_p)$ such that:
\begin{eqnarray}
\label{2.27}
	\tau(\chi) \ = \ \prod_{p \mid m} \ \tau_p(\chi_p).
\end{eqnarray}
If  $\chi',\ \chi'' \in \ZMd{m}$, then there are 
characters  $\chi_p',\ \chi_p'' \in \ZMd{p^{v(p)}}$ 
and Jacobi-Sums $j_p(\chi_p', \chi_p'')$ such that:
\begin{eqnarray}
\label{2.28}
	j(\chi', \chi'') \ = \  
		\prod_{p \mid m} \ j_p(\chi_p', \chi_p'').
\end{eqnarray}
\end{fact}
\begin{proof}
The proof is an exercise in the use of the Chinese Remainder Theorem.
\end{proof}

If $\chi$ is a primitive character, the absolute value of its Gauss-Sum 
is determined by:
\begin{eqnarray}
\label{2.30}
	\tau(\chi) \cdot \tau(\chi^{-1}) \ = \ \chi(-1) \cdot m .
\end{eqnarray}
Gauss and Jacobi-Sums are connected by:
\begin{eqnarray}
\label{2.31}
	\quad j(\chi, \chi') \cdot \tau(\chi \cdot \chi') \ = \  
		\tau(\chi) \cdot \tau(\chi'), \quad
	 \hbox{ if  $\chi$, $\chi'$ and $\chi \cdot \chi'$
		are primitive. } 
\end{eqnarray}
Let the $\chi$ be a character of conductor $m$ and order $f$; the {\sl
multiple Jacobi-Sums $J_{\nu}(\chi)$ } are defined by:
\begin{eqnarray}
\label{2.31.1}
 J_1 & = & 1 \nonumber \\ J_{\nu + 1} & = &
	J_{\nu} \cdot j(\chi, \chi^{\nu}),\quad \hbox {for $\nu = 1,
	2, \ldots, f -2$} \\ J_{f} & = & \chi(-1) \cdot m \cdot
	J_{f-1} \nonumber
\end{eqnarray}

It is easy to verify by induction that:
\begin{eqnarray}
\label{2.31.2}
 J_{\nu} & = & \frac{\tau(\chi)^{\nu}}{ \tau(\chi ^ {\nu})}, \quad \hbox {for
	$\nu = 1, 2, \ldots, f$}, 
\end{eqnarray}
where the sum of the trivial character is set by definition to $\tau(\chi^{f})=1$.
\NI The Chinese Remainder Theorem can be used for expressing the
Gauss-Sum of a character $\chi$ as a product of Gauss-Sums of
characters of prime power orders.

\begin{remark}
If $m = q$ is a prime, $\nu \in \ZMs{q}$ and $t$ and $f$ are such that
$t =\ord_q (\nu)$ and $f = \varphi(q)/t$, then $H(\nu)= \ZMs{q} \ / \,
< \nu \mod q>$ is a cyclic group isomorphic to $\{ g^{t \cdot i} \mid
i=1, 2, \ldots ,f \}$, where $g$ is a generator of $\ZMs{q}$. With
$\rho$ a primitive \nth{q} root of unity, the relations \rf{2.20} --
\rf{2.22} can be rewritten explicitly as:
\begin{eqnarray}
\label{2.34}    
  \eta_j(\rho,\nu) & = & \sum_{i=1}^{t} \ \rho^{g^{t \cdot j} \cdot \nu^{i}}
	\hbox{ for $j = 1, 2, \ldots , f$.} \\[0.3cm]
\label{2.35}    
       \tau(\chi) & = & \sum_{j=1}^{f} \ \eta_j(\rho,\nu) \cdot
		\chi(g^{t \cdot j}), \ \ \ \forall \, \chi \in
		H(\nu)^{\bot}. \\[0.3cm]
\label{2.36}    
	f \cdot \eta_j(\rho,\nu) & = & \sum_{j=1}^{f} \ \chi^{-1}(g^{t
		\cdot j}) \cdot \tau(\chi^{i}) \ \ \ \hbox{for $j = 1,
		2, \ldots , f$.}
\end{eqnarray}

It follows from \rf{2.36}, that the Gauss-Sums $\tau(\chi), \ \chi \in
 H(\nu)^{\bot}$ are {\sl Lagrange resolvents} for the Gauss Periods
 $\eta_j(\rho,\nu)$.  In this context, one can interpret \rf{2.21} as a
 generalization of Lagrange resolvents to abelian extensions. We shall
 see in the next chapter, that Gauss Periods generate intermediate
 extensions in cyclotomic fields.  The Gauss-Sums can be used to
 calculate the periods and thus to generate intermediate cyclotomic
 fields.  

Gauss sums can be defined for primitive characters of prime power
conductors; the properties arising in this context have been
investigated in \cite{Mi1} but are not of interest in our present
context. This explains the choice of $s$ as being squarefree in the
definitions above.
\end{remark}

In the case when $n = r$ is a prime and $\rg{A}$ is a field of
 characteristic $r$, the action of the Frobenius upon Gauss sums
 induces some formulae which are specific for character sums over
 finite fields. Let $\chi$ be a primitive character of conductor $m$
 and order $f$; $\zeta, \rho \in \rg{A}$ are primitive roots of unity,
 with respective orders $f$ and $m$. We investigate the action of the
 automorphism $\phi_r: x \mapsto x^r$ of \rg{A} upon $ \tau(\chi) \ :$
\[ \tau(\chi)^r \ = \ \sum_x \ \bigl(\chi(x) \cdot \rho^x \bigr)^{r} \ = \
\sum_x \ \chi^r(x) \cdot \rho^{r \cdot x}. \]
By using \rf{2.19} we have:
\begin{eqnarray}
\label{2.37}    
	\tau(\chi)^r \ = \ \chi^{-r}(r) \cdot \tau(\chi^r),
\end{eqnarray}
and iterating \rf{2.37} we get:
\begin{eqnarray}
\label{2.38}    
	\tau(\chi)^{r^k} \ = \ \chi^{-k \cdot r^k}(r) \cdot \tau(\chi^{r^k})
	\hbox{, for $k \geq 1$ .} 
\end{eqnarray}
If $r^t = 1 \mod f$, then
\begin{eqnarray}
\label{2.39}    
	\tau(\chi)^{r^t - 1} \ = \ \chi^{-t}(r).
\end{eqnarray}
The relations \rf{2.37} and \rf{2.39} are central in primality
testing.  It will be important to have efficient computing methods for
powers of Gauss and Jacobi sums, if they are to be used in practical
algorithms.

\section{Further Criteria for Existence of Cyclotomic Extensions}
 
The condition $(II)$ in Theorem \ref{th1} is central for primality
 proving and motivates the interest in proving the existence of
 cyclotomic extensions. One way of doing this is shown in Theorem
 \ref{LL} and it generalizes the classical Lucas -- Lehmer tests. The
 condition $(III)$ can be connected by relation \rf{2.22} to Gauss
 periods and sums.

The resulting conditions indicate the direction for the Jacobi sum
test. Before stating them, let us introduce some notations. Let $n, s,
t$ be like in Theorem \ref{th1} and $\xi_t, \xi_s \in \C$ be fixed;
furthermore, we assume that there exists a saturated \nth{t}
cyclotomic extension $\rg{R} \supset \id{N}$ and $\zeta \in \rg{R}$ is
a primitive \nth{t} root of unity. We shall write like previously
$\ZMd{s}$ for the characters with image in $\rg{R}$ while $ \ZMb{s} =
\left\{ \ \chi : \ZMs{s} \rightarrow < \xi_t > , \hbox{ with $\chi$
multiplicative} \ \right\}$. For $a \in \ZMs{s}$ we let $H(a) =
\ZMs{s}/< a \bmod s>$ and
\[ H(a)^{\widehat{} \ } = \left\{ \chi \in \ZMd{s} \ : \ \chi(a) = 1 \right\}
\subset \ZMd{s} \] be its dual. The set $H(a)^{\top} \subset \ZMb{s}$ is
defined by analogy. Then,

\begin{theorem}  
\label{th2}
The following statement is equivalent to $(I) - (IV)$ of Theorem \ref{th1}:
\begin{itemize}
\item[(V)] If the Gauss sums $\tau(\chi)$ are defined for $\chi \in \ZMb{s}$
            with respect to $\xi_s$, then: 
\[ \chi \in H(n)^{\bot \ } \ \Longleftrightarrow \ \exists \ \hbox{ a
            homomorphism } \quad \vartheta: \Z[\xi_t,\tau(\chi)]
            \rightarrow \rg{R}.\]
\end{itemize}
\end{theorem}
\begin{proof}
Suppose that $(III)$ holds, thus a map $\tau_0: \rg{A} =
\id{O}\left(\Q(\xi_s)^{< n \bmod s >}\right) \rightarrow \id{N}$
exists. In particular, it follows that the Gauss periods
$\eta_h(\xi_s, n) = \sum_{\mu \in h} \ \xi_s^{\mu}$ with $h \in H(n)$
are mapped to $\id{N}$. Let $\vartheta$ be the lift of $\tau_0$ with
$\vartheta(\xi_t) = \zeta \in \rg{R}$. If $\tau(\chi)$ are Gauss sums
with respect to $\xi_s$ and $\chi \in H(n)^{\bot}$, then we gather
from \rf{2.21} that $\vartheta\left(\tau(\chi)\right) \in \rg{R}$,
which proves that $(III) \Rightarrow (V)$.

Suppose now that $(V)$ holds and let $\rg{B} \subset \Z[\xi_t, \xi_s]$
be the ring generated by $\xi_t$ and the Gauss sums $\tau(\chi), \
\chi \in H(n)^{\bot}$, while $\vartheta : \rg{B} \rightarrow \rg{R}$
is such that $\vartheta(\tau(\chi)) \in \rg{R}$. Using \rf{2.22} we
see that $\vartheta$ maps the Gauss periods $\eta_h$ to $\rg{R}$, and
if $\sigma$ generates the Galois group of $\rg{R}/\id{N}$ acting on
$\zeta$, then $\vartheta(\eta_h)$ are $\sigma$ invariant, so
$\vartheta(\eta_h) \in \id{N}$. Using reduction modulo primes $r | n$
and arguments from the proof of Theorem \ref{th1}, we deduce that $r
\in < n \bmod s >$ and thus $(V) \Rightarrow (II)$, which completes
the proof.
\end{proof}
Note that since only characters $\chi \in H(n)^{\bot}$ are considered,
the condition $(V)$ is a slight improvement of the one used in the
initial form of the Jacobi sum test \cite{APR}, and which involved all
characters in $\ZMd{s}$.
\begin{lemma}
\label{lcole}
Let $p, q$ be primes not dividing $n$, with $p^k \parallel (q-1)$ and
$(\rg{R}, \sigma, \zeta)$ be a saturated \nth{p} cyclotomic extension
of $\id{N}$. Let $\ \chi \in \ZMd{q}$ be a character of order $p^k$
and $\alpha, \beta \in \rg{R}$ be given by:
\begin{eqnarray}
\label{4.14}
 \alpha & = & J_{p^k}(\chi) \quad \hbox{ and} \\ \nonumber \beta &
	 = & J_{\nu}(\chi), \quad \hbox{ where $\nu = n \mod
	 p^k$.}  \nonumber
\end{eqnarray}
Let $l = [ n / p^k]$ and suppose that
\begin{eqnarray}
\label{4.15}
   \alpha^l \cdot \beta = \eta^{-n} \quad \hbox{holds for some $\eta
	\in <\ \zeta \ >$. }
\end{eqnarray}
Then $\eta=\chi(n)$ and $\chi(r)=\chi(n)^{l_p(r)}, \  \forall r \mid n$ ,  with
$l_p(r)$ defined in Lemma \ref{lp}.
\end{lemma}
\begin{proof}
Let $\rg{R}' = \rg{R}[X]/(\Phi_q(X))$ and define $\zeta_q = X +
\Phi_q(X) \in \rg{R}'= \rg{R}[X]/(\Phi_q(X))$; one proves that
$\rg{R}'$ is a galois extension of $\rg{R}$ and also of $\id{N}$. We
then define the Gauss sum $\tau(\chi)$ with respect to $\zeta_q$ and
claim that the identities on multiple Jacobi sums hold for this sum;
this is a simple verification and is left to the reader. The actual
identities are meaningful in the ring $\rg{R}$, but we need $\rg{R}'$
for introducing the Gauss sums. By the definition of $\alpha, \beta$
and $l$, \rf{4.15} is equivalent to
\begin{eqnarray}
\label{4.16}
	\tau( \chi )^n \ = \ \eta^{-n} \cdot \sigma(\tau( \chi ) ).
\end{eqnarray}
Raising \rf{4.16} to the power $n$ repeatedly, we find:
\begin{eqnarray}
\label{4.17}
\tau( \chi )^{n^i} \ = \ \eta^{-i \cdot n^i}
	\cdot \sigma^i (\tau( \chi ) ) \quad \forall i \geq 1
\end{eqnarray}
and, with $i = p^k \cdot ( p-1)$ and $N = n^i$ ,
\begin{eqnarray}
\label{4.18}
	\tau( \chi )^{N-1} = 1.
\end{eqnarray}
If $r \mid n$ is a prime and $\eu{R} \subset \rg{R}'$ a maximal ideal
through $r$, then by \rf{lpr}
\begin{eqnarray}
\label{4.19}
	\tau( \chi )^{r}\ = \ \chi(r)^{-r} \cdot (\tau( \chi^r ) )
		\mod \eu{R}'.
\end{eqnarray}
From the existence of the saturated \nth{p} extension $\rg{R}$ we
gather, by Fact \ref{lp}, that there are two integers $l_p(r), u_p(r)$
verifying \rf{lpr}. With these, we let $m \in \N$ be such that
$m=l_p(r) \mod p^k$ and $m=u_p(r) \mod (p-1)$, so that $\sigma^m( \chi
) = \chi^r$ and
\begin{eqnarray}
\label{4.20}
	v_p(r-n^m)\  = \  v_p \bigl(n^m \cdot (r / n^m -1) \bigr) 
		\geq v_p(N-1)	
\end{eqnarray}
We let $i = m$ in \rf{4.17}, use $\sigma^m \bigl(\tau(\chi)\bigr) \ =
\ \tau(\chi^r)$ and divide by \rf{4.19}. This is allowed, since $\tau(
\chi ) \cdot \tau(\chi^{-1}) = \pm q $ and $( q, n ) = 1$; the result
is:
\begin{eqnarray}
\label{4.21}
	\tau( \chi )^{n^m-r} \ = \ \bigl( \chi(r) \cdot \eta^{-m}
		\bigr)^r \bmod \eu{R}.
\end{eqnarray}
Let $u$ be the largest divisor of $(N-1)$ which is coprime to
$p$. From \rf{4.18}, \rf{4.20} and by raising \rf{4.21} to the power
$u$, we get:
\begin{eqnarray}
\label{4.22}
 1 \ = \ \bigl( \chi(r) \cdot \eta^{-m}
	\bigr)^{r \cdot u} \mod \eu{R}.
\end{eqnarray} 
Now $\rho = \chi(r) \cdot \eta^{-m} \in \rg{R}$ is a primitive root of
unity of some order $p^v$ and such that $\rho \equiv 1 \bmod
\eu{R}$. We claim that $v = 0$ and $\rho = 1$; if this was not the
case, then $p^v = \prod_{i=1}^{p^v - 1} (1 - \rho^i) = \frac{X^{p^v} -
1}{X - 1} \bigg|_{X = 1}$ and since $\rho \equiv 1 \bmod \eu{R}$, we
should have a fortiori $p^v \equiv 0 \bmod \eu{R}$ which contradicts
$(p, r) = 1$. So $\rho = 1$ and thus $\chi(r) = \eta^{m} =
\eta^{l_p(r)}$.  This holds for all primes $r \mid n$ and, by
multiplicativity, for all divisors $r' \mid n$. In particular, since
$l_p (n) = 1$, it follows that $\eta = \chi( n )$.
\end{proof}
\begin{remark}
The equivalent relations \rf{4.15} and \rf{4.16} are reminiscent of the
identity \rf{2.37} holding in finite fields. The statement of the Lemma holds
a fortiori when replacing \rf{4.15} by
\begin{eqnarray}
\label{fullpow}
\alpha^{(n^{t_p}-1)/p^k} = \chi^{-t_p}(n), \quad \hbox{ with } \quad t_p =
\ord_{p^k}(n),
\end{eqnarray}
which is the analog of \rf{2.39} and is obtained by iteration of \rf{4.16}.
Here $\alpha = J_{p^k}(\chi)$ like in the hypothesis above.
\end{remark}
The Lemma \ref{lcole} indicates the steps for proving the existence of
\nth{s} cyclotomic extensions with Jacobi sums. This is the corner stone
of the Jacobi sum test:
\begin{corollary}
\label{jstest}
Suppose that $s$ is square-free, $t = \ord_s(n)$ and $\rg{R}$ is a
saturated \nth{t} extension of $\id{N}$ with $\zeta \in \rg{R},
\Phi_t(\zeta) = 0$. We let $\ZMd{s}$ be the set of characters of
conductor $s$ with images in $< \zeta >$, the sets $\id{P}, \id{Q}$ be
given by \rf{sets} and
\begin{eqnarray}
\label{chars}
\quad  \id{C} = \left\{ \ \chi_{\wp} \in \ZMd{s} \ : \ \wp \in \id{Q}, \
\chi \ \hbox{ has conductor $q$ and order $p^k$ } \right\}.
\end{eqnarray}
Suppose that
\begin{eqnarray}
\label{jt}
\tau(\chi_{\wp})^{n - \sigma} \in \ < \chi_{\wp}(n) >, \quad \forall \
\wp \in \id{Q}.
\end{eqnarray}
or, alternately, for all $\wp \in \id{Q}$ one has:
\begin{eqnarray}
\label{jt1}
\quad \alpha_{\wp}^{\frac{n^{t_p} - 1}{p^k}} =
\chi_{\wp}(n)^{-t_{\wp}}, \quad \hbox{ with } \quad t_{\wp} =
\ord_{p^k}(n) \quad \hbox{ and } \quad \alpha_{\wp} =
\tau(\chi_{\wp})^{p^k}.
\end{eqnarray}
Then an \nth{s} cyclotomic extension of $\id{N}$ exists.
\end{corollary}
\begin{proof}
Using Lemma \ref{lcole}, respectively \rf{fullpow}, we deduce from \rf{jt} or
\rf{jt1} that $\chi(r) = \chi\left(n^{l_p(r)}\right)$ for all the characters
$\chi \in \ZMd{s}$. Let $L(r) \equiv l_p(r) \bmod p^k$ for all $p^k \parallel
t$; then we have a fortiori $\chi(r) = \chi\left(n^{L(r)}\right)$ for all
$\chi \in \ZMd{s}$ and by duality, $r \equiv n^{L(r)} \bmod s$. This holds for
all $r | n$ which implies $(II)$ and the fact that an \nth{s} cyclotomic
extension of $\id{N}$ exists.
\end{proof}

The conditions for existence of \nth{s} cyclotomic extensions, which are based
on Gauss sums, require $s$ to be squarefree. This is not the case for the
Lucas -- Lehmer test in Theorem \ref{LL}. We wish to combine the information
about extensions proved by the two methods. This happens to be quite easy,
since the extensions proved by means of Theorem \ref{LL} are saturated and
thus \rf{lpr} holds by Lemma \ref{lp}. We group these observations in 
\begin{fact}
\label{comb}
Let $(s_1, s_2) = 1$ with $s_2$ squarefree, $s = s_1 \cdot s_2$ and $t_i =
\ord_{s_i}(n), i = 1, 2$, $t = \ord_s(n)$. Suppose that $(\rg{R}, \sigma,
\zeta)$ is a saturated \nth{t} cyclotomic extension of $\id{N}$ and $t_1 | [
\rg{R} : \id{N} ]$. Furthermore there is a $\beta \in \rg{R}$ with
\[  \Phi_{s_1}(\beta)  =  0, \quad \hbox{ and } \quad \beta^n  =  
\sigma(\beta),\] 
such that $(\rg{R}, \beta, \sigma)$ is saturated as a \nth{s_1} extension. If
the conditions of Corollary \ref{jstest} apply for $s = s_2$, then an \nth{s}
cyclotomic extension of $\id{N}$ exists.

Furthermore, if $s_i$ are any coprime integers such that saturated
\nth{s_i} cyclotomic extensions of $\id{N}$ exist and $s = \prod_i
s_i$, then a saturated \nth{s} extension exists.
\end{fact}
\begin{proof}
Let $r | n$ be a prime. The proof of Corollary \ref{jstest} and the
fact that the \nth{s_1} extension is saturated imply, by means of
Lemma \ref{lp}, that $\chi(r) = \chi\left(n^{L(r)}\right)$ for all
characters $\chi \in \ZMd{s}$ with $L(r) \equiv l_p(r) \bmod
p^{v_p(t)}$ and all $p | t$. The statement about combinations of
saturated extensions is a direct consequence of Lemma \ref{lp}.
\end{proof}

\section{Certification}
Certificates for primality proofs are data collected during the
performance of the test of primality for a given number $n$. The
certificate allows to perform a verification of the primality of $n$
in (sensibly) less time than it took to collect the data. A recursive
{\em Pratt} certificate \cite{Pr} is the following: suppose that $n =
a F + 1$ is a prime and $\prod_{i=1}^k p_i^{e_i} = F > \sqrt{n}$, with
$q_i = p_i^{e_i}$ being prime powers. Furthermore, suppose that $b_i
\in \Z$ are such that $\Phi_{q_i}(b_i) \equiv 0 \bmod n$, or $b_i
\equiv c_i^{(n-1)/q_i} \bmod n$, while $(c_i^{(n-1)/p_i}-1,n ) = 1$
and $c_i^{(n-1)/p_i} \equiv 1 \bmod n$. A certificate $C(n)$ is
defined recursively by
\[ \ C( n ) = \{ \ b_i : i = 1, 2, \ldots, k \ \} \ \bigcup \ 
\left( \cup_{i=1}^{k} \ C(p_i) \right) , \] with $C( p_i ) $ being
certificates for $p_i$. If $b_i$ are computed by trial and error,
using the $c_i$ above, the time for building a certificate is larger
than the one required for its verification.  This example suggests a
generalization to CPP. We may mention that it was
believed until recently that certification was an advantage of ECPP
and not achievable for CPP. It is not the case, as we show here.

The relation \rf{jt1} shows that if $\chi_{\wp}(n) = 1$ and an \nth{s}
cyclotomic extension does exist, then not only $\tau(\chi_{\wp}) \in
\rg{R}$ as follows from $(IV)$, but it also can be explicitely
computed in $\rg{R}$ by means of the Theorem \ref{troot}. This would
provide for a certificate which can be verified by exponentiations
with exponent $p^k$ in $\rg{R}$; however the list $\id{C}$ contains
also characters which do not vanish at $n$. In such cases, one first
modifies $\alpha_{\wp}$ accordingly before taking a \nth{p^k} root. 

The resulting criteria are given in 
\begin{theorem}
\label{cert}
Let $s$ be squarefree, $t = \ord_s(n)$ and $\rg{R} = \id{N}[\zeta]$ be
a saturated \nth{t} cyclotomic extension. Let $\id{Q}, \id{P}, \id{C}$
be defined in \rf{sets}, \rf{chars} and suppose that for all $\wp \in
\id{Q}$ there is a $\beta_{\wp} \in \rg{R}$ such that, 
for\footnote{We supress here, for typographic reasons, writing out the explicite dependency on $\wp$.} $t = t_{\wp} = \ord_{p^k}(n)$:
\begin{eqnarray}
\label{cercon}
\quad \quad \beta_{\wp}^{p^k} = \chi_{\wp}(n)^{\frac{t \cdot
p^k}{n^{t} -1}} \cdot \alpha_{\wp}, \quad \hbox{ with } \quad
t = \ord_{p^k}(n) \quad \hbox{ and } \quad \alpha_{\wp} =
\tau(\chi_{\wp})^{p^k}.
\end{eqnarray}
Then an \nth{s} cyclotomic extension of $\id{N}$ exists.
\end{theorem}
\begin{proof}
If $n$ is prime, then
\[  \left( \alpha_{\wp} \cdot \chi_{\wp}(n)^{\frac{t \cdot
p^k}{n^{t} -1}} \right)^{\frac{{n^{t} -1}}{p^k}} = 1 \] as
a consequence of \rf{2.39} and the expression $
\chi_{\wp}(n)^{\frac{t \cdot p^k}{n^{t} -1}} \cdot
\alpha_{\wp}$ is in this case a \nth{p^k} power in $\rg{R}$, as
follows from Theorem \ref{troot}. The existence of $\beta_{\wp}$ is a
necessary condition for primality and thus consistent with our
purpose.

Since $\rg{R}$ is saturated, $\rg{S} = \rg{R}[X]/\left(X^t -
\zeta\right)$ is a \nth{t^2} cyclotomic extension and in particular
galois with group of order $t \cdot \ord_t(n)$. We claim that $\rg{S}$
contains a primitive \nth{s} root of unity $\omega$ upon which
$\vartheta$ acts making $(\rg{S}, \omega, \vartheta)$ into an \nth{s}
cyclotomic extension in the sense of Remark \ref{subext}. Our proof
relays upon Theorem \ref{th2}.

We first prove an auxiliary fact about saturation. Let $\wp \in
\id{Q}$, let $$\delta_{\wp} =
\chi_{\wp}(n)^{-\frac{t}{n^{t}-1}}$$ and $u = k_p(n)$ be
the saturation exponent of $p$ with respect to $n$. Then there is an
integer $0 \leq v < p^u$ such that
\[   \frac{t}{n^{t}-1} = \frac{v}{p^u} + m, \quad \hbox{ with } 
\quad m \in \Z, \]
and hence 
\begin{eqnarray}
\label{nosat}
\delta_{\wp} \times \chi_{\wp}(n)^{\frac{v}{p^u}} \in \rg{R} \quad
\hbox{ and } \quad \delta_{\wp} \in \rg{S}.
\end{eqnarray}
Indeed, let $t' = \ord_p(n)$ so that $v_p(n^{t'} - 1) = u$ and
 suppose $\frac{t'}{n^{t'}-1} = \frac{v}{p^u} \bmod \Z$. The
 assertion follows for $\wp = ( p^k , q)$ with $k \leq u$; if $k = u +
 j$, then by the definition of saturation, $t =
 \ord_{p^{u+j}}(n) = p^j \cdot t'$. Since
\[  \frac{n^{p^j \cdot t'} - 1}{p^j (n^{t'} - 1) } \equiv 1 \bmod p^u,  \] 
as shown by a short calculation, it follows that
\[ \frac{t}{n^{t}-1} = \frac{t'}{n^{t'}-1} \times 
\left(p^j \cdot \frac{n^{t'} - 1}{n^{p^j t'} - 1}\right) \equiv
\frac{v}{p^u} \bmod \Z , \] thus proving the claim. Note that
$\gamma_{\wp} = \beta_{\wp} \cdot \delta_{\wp}$ is a solution of
$X^{p^k} = \alpha_{\wp} = \tau\left(\chi_{\wp}\right)^{p^k}$. 

Let as usual $\xi_t, \xi_s \in \C$ be fixed and $\psi \in H^{\top}(n)$
be a character of order $p^k$ with image in $< \xi_t >$, satisfying
$\psi(n) = 1$ and let $\chi \in \ZMd{s}$ be the image of $\psi$ by
$\theta : \xi_t \mapsto \zeta$. We want to show that $\theta$ can be
extended to $\tau(\psi)$.  This is done as follows: $a(\chi) =
\left(\tau(\psi)\right)^{p^k} \in \Z[\xi_t]$ so we can set
$\alpha(\chi) = \theta(a(\psi)) \in \rg{R}$ and then $\Z[\xi_t,
\tau(\psi)] \subseteq \Z[\xi_t, X]/( X^{p^k} - a(\psi) )$. The map
$\theta$ extends to $\tau(\psi)$ if we can show that the equation
$T^{p^k} = \theta(a(\chi)) = \alpha(\chi)$ has a solution in
$\rg{R}$. Furthermore, if this holds for any $\wp \in \id{Q}$, we
conclude that for each $\psi \in H^{\top}(n)$, the Gauss sum
$\tau(\psi)$ maps to $\rg{R}$ and the claim then follows from
$(V)$. Now if $\psi \in H^{\top}(n)$ is a character of order $m$, it
can be decomposed in a product of characters $\psi = \prod_{p^k
\parallel m} \ \psi_p$ of characters of prime power orders $p^k ||
m$. The Gauss sum $\tau(\psi) = J(\psi) \times \prod_{p | m}
\tau(\psi_m)$ where we assumed that $\theta \left( \tau(\psi_m)
\right) \in \rg{R} $ and $J(\psi)$ is a product of Jacobi sums which
also maps to $\rg{R}$.

Suppose that the prime decomposition of $s$ is $s = \prod_{i=1}^u q_i$
and define the factor characters $\chi_i(x) = \chi(x \bmod q_i)$; the
decomposition formula \rf{2.27} implies that $\tau(\chi) =
\prod_{i=1}^d \ \tau(\chi_i)$. By definition of $\id{P}$, there are
pairs $\wp_i = (p^{k_i}, q_i) \in \id{P}$ such that $\chi_i =
\chi_{\wp_i}$. Using \rf{2.19}, we have
$\left(\tau(\chi_i)\right)^{p^{k_i}} = \left(\beta_i \cdot
\delta_i\right)^{p^{k_i}} \in \rg{R}$ and $\beta_i =
\beta_{\wp_i}$, etc. Note that we have to raise to the power $p^{k_i}$
in the previous formula, in order to consider elements which are
defined in $\rg{R}$; an alternative solution would be a formal
adjunction of an \nth{s} root of unity to $\rg{R}$. The hypothesis
$\chi(n) = 1$ and relation \rf{nosat} imply that
\[ \tau(\chi)^{p^u} = \left( \prod_{i=1}^d \beta_i \right)^{p^u} \times 
 \prod_{i = 1}^d \ \chi_i(n)^{v + m_i p^u} = \left(\prod_{i=1}^d
\beta_i \right)^{p^u} \cdot \chi(n)^{v} \cdot \prod_{i=1}^d
\chi_i(n)^{ m_i p^u} = \beta^{p^u} , \] with $m_i = \frac{t_i}{n^{t_i}
- 1} - \frac{v}{p^u}$ and $\beta = \prod_{i=1}^d \beta_i \cdot
\chi_i(n)^{m_i} \in \rg{R}$. This shows that
$\vartheta\left(\tau(\psi)\right) \in \rg{R}$ as claimed, and
completes the proof for odd $p$ or $p = 2$ and $n \equiv 1 \bmod
4$. If $n \equiv 3 \bmod 4$ and $p = 2$, the saturation context is
different. The proof uses an appropriate variant of \rf{nosat} and 
shall be skipped here.
\end{proof}

It is useful to note, that \rf{fastpow} substantially accelerates the
evaluation of \rf{jt1}, making it comparable to the one of \rf{jt}. As
a consequence, computing a certificate requires no substantial
additional work compared to the classical Jacobi sum test.

\subsection{Computation of Jacobi Sums and their Certification}
We are interested in the computation of Jacobi sums $j(\chi, \chi^a)$,
where $\chi = \chi_{\wp}$ is a character of prime conductor $q$ and prime
power order $p^k | (q-1)$. For these sums, the absolute value is
\begin{eqnarray}
\label{normjac}
 j(\chi, \chi^a) \times \overline {j(\chi, \chi^a)} = q .
\end{eqnarray}

Since the conductor $q$ of Jacobi sums in CPP has superpolynomial size,
their computation is a critical step which deserves some attention. From 
the theoretical point of view, the recent random polynomial algorithm of 
Ajtai, Kumar and D. Sivakumar \cite{AjKuSi}  for 
finding shortest vectors in lattices solves the concrete problem 
in polynomial, and in fact linear time and space. Indeed, as we detail below, 
Jacobi sums of characters of order $P$ are shortest vectors in certain {\em well rounded 
lattices}, i.e. lattices with a base of vectors of equal length. In a lattice of dimension 
$P$, the algorithm \cite{AjKuSi} takes $O(2^P)$ space and time, and since in the context
of CPP, the size $P = O(\log_{(2)}(P)$, it follows that Jacobi sums can be computed 
in random linear time. 

In practice, the dimensions of lattices are quite small and in view both of constants and
implementation complexity of the shortest vector algorithm, it is useful to discuss
some simpler practical methods too.

For moderate values of $q$, possibly $q < 10^{14}$, the direct
computation based on the definition \rf{2.25} is adequate and
fast. The bottleneck is the necessity to store a table of discrete
logarithms modulo $q$. This can simply be avoided, by performing the 
computation of Gauss periods in $\C$, then computing Gauss and Jacobi sums
in $\C$ too; finally, from the conjugates of a Jacobi sum, one recovers 
its coefficients as an algebraic integer. The method is straightforward and was 
implemented in the Master Thesis \cite{Ma}.

For larger conductors, it is preferable to use methods of
lattice reduction. These have been investigated in \cite{BuKo},
\cite{Mi1}, \cite{Wa} and are based on the following observation. Let
$\eu{Q} \subset \Z\left[\xi_{p^k}\right]$ be a prime ideal above $q$;
note that the choice of $p$ implies that $q$ splits completely and
$\eu{Q}$ has inertial degree one. Let $G =
\Gal\left(\Q(\xi_{p^k})/\Q\right)$ and $I = \Z[ G ]$ be the
Stickelberger ideal. There is an element
\begin{eqnarray}
\label{stick}
 \theta = \sum_{(c, p) = 1; \ 0 < c < p^k } \left[\frac{ac}{p^k} \right] 
\cdot \sigma_{c}^{-1}\in I 
\end{eqnarray}
such that
\[  \left( j(\chi, \chi^a) \right) = \eu{Q}^{\sigma \theta}, \]
for some $\sigma \in G$. The ideal $\eu{Q}^{\theta}$ can be
represented by a $\Z$ - base, being a free $\Z$ - module of rank
$\varphi(p^k)$. As such, it is a lattice and it follows from
\rf{normjac} that $\sigma^{-1} \left(j( \chi, \chi^a )\right) \in
\eu{Q}^{\theta}$ is a shortest vector of this lattice, with respect to
the embedding (Gauss) norm $\parallel x \parallel \ = \ \sum_{\sigma
\in G} | \sigma(x) |^2$. 

This opens the road for applications of
methods of lattice reduction. Without entering in details, which can
be found in the references, we mention that lattice reduction allows
use of large conductors, but the growth of the order -- which controls
the dimension of the lattice -- is critical. Indeed, the problem of
finding the shortest vector in a lattice of dimension $d$ with initial
base of vectors bounded by $q$ has complexity $O\left( d^d \cdot \log
(q)^{O(1)} \right)$. In practice, due in part to the particularity
that the lattices to consider (generated by Jacobi sums) have a basis
of shortest vectors -- they are {\em well rounded} -- the computations
are quite efficient, and shortest vectors are frequently found directly 
by LLL, for character orders up to at least $P \sim 125$ \cite{Ma}. 

A more efficient LLL based approach which works for small class numbers of the 
cyclotomic field $\Q(\zeta_P)$ follows the method used by Buhler and Koblitz 
in \cite{BuKo}: Let $\eu{Q} \subset \Q(\zeta_P)$ be an ideal above the conductor
$q$. If $h$ is the class number, then find by LLL a generator of $\eu{Q}^h$ 
and compute Jacobi sum powers $j(\chi, \chi')^h$ by use of Stickelberger elements. 
If the generator of $\eu{Q}^h$ is found correctly by LLL, then this method uses
only one LLL computation for a given conductor and order.

Finally, the implementations of PARI for computing the structure of class and unit 
groups of number fields turned out to be very efficient in computing Jacobi sums too. The bottle neck there is the space requirement, since finding generators of principal ideals
is based on building up all the information on class and unit groups. Here, we use the fact that multiple Jacobi sums have to be computed in the same field, so the field construction which is slower, happens only once.

Since the computations in $\C$ and the LLL based method are not guaranteed to yield Jacobi sums -- the first due to rounding errors, the second due to the shortest vector problem -- 
it is therefor interesting that one can {\em certify} very easily that
the value of a Jacobi sum is correct, using the very formulae
displayed above. This comes both as a verification and as part of a ceritificate for ulterior verifications of a primality proof. More precisely, we have:
\begin{lemma}
\label{certjac}
Let $\wp = (p^k, q)$ with $p^k | (q-1)$ and $p, q$ being primes. Let
$0 < a < p^k$ be an integer and $\xi = \xi_{p^k} \in \C$ be fixed. If
$\alpha \in \Z[\xi]$, there is a deterministic algorithm which
verifies whether $\alpha = j(\chi, \chi^a)$ for some character $\chi
\in \ZMb{q}$ of conductor $q$ and order $p^k$. The verification is
done in $\id{O}\left(p^{2k} \cdot \log(q)\right)$ binary operations.
\end{lemma}
\begin{proof}
A first condition which must be fulfilled by a Jacobi sum is the local
$p$ - adic norming condition $j \equiv \pm 1 \bmod (1 - \xi)^2$, see
e.g. \cite{IR}, and this fixes the choice of a root of unity factor
\footnote{The sign is always positive, if one adopts Lang's definition 
of Gauss sum, with a minus sign.}. Thus one starts by verifying that
\begin{eqnarray}
\label{initc}
\alpha \times \overline{\alpha} = q, \quad \hbox{ and } \quad \alpha
\equiv \pm 1 \bmod (1-\xi)^2,
\end{eqnarray}
in $O(p^k \cdot \log(q))$ operations -- note that the coefficients of
a Jacobi sums have size $\sim \ \sqrt{q}$ and thus a multiplication of
two Jacobi sums has the complexity above.

Since $q \equiv 1 \bmod p^k$ there is a $c \in \Z$ with $\Phi_{p^k}(c)
\equiv 0 \bmod q$ and thus $\eu{Q} = (\xi - c, q)$ is a prime ideal
above $q$. Next one computes $\beta = (\xi-c)^{\theta} \in \Z[\xi]$
with $\theta \in I$ defined by \rf{stick}. This is done in $O(p^{2k}
\cdot \log(q))$ operations. Finally, one checks if there is a $\sigma
\in G$ such that $\sigma(\beta) \equiv 0 \bmod \alpha$. If yes, then
$\alpha$ is a Jacobi sum and $\alpha = j(\chi, \chi^a)$ for some
character of order $p^k$ and conductor $q$, otherwise the claim is
false.
\end{proof}

\section{Algorithms}
The previous sections provide the theoretical foundation for the CPP
primality proving algorithms. These consist of three steps, which are
partially interdependent. Like usual, we denote by $n$ a number to be
proved prime and $\id{Q}, \id{P}, \id{C}$ are defined by \rf{sets} and
\rf{chars}, respectively. The main steps of the algorithms are the
following:
\begin{itemize}
\item[A.] \textbf{Work Extensions:} Select two parameters $s, t$ such
that $t = \ord_s(n)$ and build a saturated \nth{t} extension
$\rg{R}/\id{N}$ -- e.g. by using the Lucas -- Lehmer method of Theorem
\ref{LL}.
\item[B.] \textbf{Parameters:} Let $s' | (n^t - 1 )$ be a totally
  factored part\footnote{It is assumed that $s'$ is built up from
  primes $q' | s'$ such that the orders $t(q') = \ord_{q'}(n) | t$ are
  \textit{small}} with $(s, s') = 1$, let $s_1$ be the order of a
  saturated \nth{s'} extension -- thus $s_1 = \prod_{q' | s' } \
  q^{k_{q'}(n)}$ and $S = s \cdot s_1$. Verify $S > \sqrt{n}$. An
  optimization cycle can lead back to A. At the end, optimal
  values of $S, s, s'$ and $t$ are chosen and the fixed conditions $S
  > \sqrt{n}$ and $t = \ord_s(n)$ hold.
\item[C.] \textbf{Test part:} 
\begin{itemize}
\item[C1.] Prove the existence of a saturated \nth{s'} cyclotomic
extension in $\rg{R}$, by using the Theorem \ref{LL}. This is the
Lucas -- Lehmer part of the test, and it can be void.
\item[C2.] Build the sets $\id{P}, \id{C}$, with respect to the
current value of $s$ and verify \rf{4.15} for all characters $\chi \in
\id{C}$. This is the Jacobi sum part of the test.
\item[C2'.] Alternately, if a certificate is required along with the
test, after building the list $\id{C}$, one finds $\beta_{\wp} \in
\rg{R}$ verifying \rf{cercon}.
\item[C3.] Perform the final trial division, verifying that
  \rf{findiv} yields no nontrivial factors of $n$.
\end{itemize}
\end{itemize}
Unless $n$ has some special form, so that many prime factors of $F_k =
n^k - 1$ are known for small $k | t$, the parameter $s'$ is either set
to $1$ and thus neglected, or gained by investing some time in the
factorization of the same $F_k$. An important observation, which does
not influence the asymptotic behavior of the algorithms but generates
a useful speed up, consists of the fact that one can verify \rf{4.15}
\textit{simultaneously} for a set of characters of mutually coprime
orders.

\begin{definition}
We define an \textit{amalgam} as a subset $\id{A} \subset \id{Q}$ such
that $\{ p(\wp) \ : \wp \in \id{A} \}$ are pairwise coprime. If $\wp =
(p^k, q)$ and $t(\wp) = \ord_{p^k}(n)$, then an amalgam $\id{A}$ is
\textit{rooted}, if there is a $\wp_0 \in \id{A}$ such that $t(\wp) |
t(\wp_0)$ for all $\wp \in \id{A}$.
\end{definition}
The relevance of amalgams is provided by the following:
\begin{theorem}
\label{thamalgam}
Let $\id{A}$ be an \textit{amalgam}, 
\[ f = f(\id{A}) = \prod_{\wp \in
\id{A}} \ p^k(\wp), \quad f'=\rad{f} = \prod_{\wp \in \id{A}} \
p(\wp), \quad t = \ord_{f'}(n) , \] and $(\rg{R}, \sigma, \zeta)$ a
saturated \nth{f} cyclotomic extension of $\id{N}$, the roots $\{
\zeta_p : p = p(\wp), \wp \in \id{A}\} \subset \rg{R}$ being all
saturated of orders $p(\wp)$. For $\wp \in \id{A}$, let:
\begin{eqnarray}
\label{4.25}
   \alpha(\wp) & = & J_{p^{k}(\wp)}(\chi_{\wp}) \quad \hbox{ and} \\
		\nonumber \beta(\wp) & = & J_{\nu(\wp)}(\chi_{\wp}),
		\quad \hbox{ where } \ \nu(\wp) = n \ \rem \ p^k.
\end{eqnarray}
Let $n = f \cdot l + \nu$ with $0 \leq \nu < f$ and $\nu = p^k(\wp) \cdot
\lambda(\wp) + \nu(\wp)$, for $\wp \in \id{A}$. Define $\alpha$ and $\beta$ by
\begin{eqnarray}
\label{albet}
   \alpha & = & \prod_{\wp \in \id{A}} \ \alpha(\wp)^{f/p^k(\wp)} \in \rg{R}
		\quad \hbox{ and} \\ \nonumber \beta & = & \prod_{\wp \in
		\id{A}} \ \alpha(\wp)^{\lambda(\wp)} \cdot \beta(\wp) \in
		\rg{R}.
\end{eqnarray}
Suppose there is an $\eta \in \ < \zeta_f >$ such that
\begin{eqnarray}
\label{4.27}
   \alpha^l \cdot \beta \ = \ \eta^{-n}
\end{eqnarray}
Then $\chi(\wp) ( r ) = \chi(\wp) (n )^{l_p(r)}, \ \ \forall r \mid n$
and $\wp \in \id{A}$. Furthermore $\eta = \prod_{\wp \in \id{A}} \
\chi(\wp)( n )$.
\end{theorem}
\begin{proof}
The proof is similar to the one of the Lemma \ref{lcole}. We shall
describe the general ideas and refer the reader to \cite{Mi1} for the
complete proof. One first adds formal \nth{q(\wp)} roots of unity to
$\rg{R}$ in order to define some Gauss sums which verify, by
definition of $\alpha$ and $\beta$ and \rf{4.27}:
\begin{eqnarray}
\label{4.28}
   	\prod_{\wp} \ \tau(\chi_{\wp})^n = \eta(\wp)^{-n} \cdot \prod \
   	\tau(\sigma(\chi_{\wp})).
\end{eqnarray}
Then one decomposes $\eta$ in a product of \nth{p} power roots of
 unity and raising \rf{4.28} repeatedly to the \nth{n} power, obtains:
\begin{eqnarray}
\label{4.29}
   \prod_{\wp} \ \tau(\chi_{\wp})^{n^h}\ = \prod_{\wp} \ \eta(\wp)^{-h \cdot
	n^h} \cdot \sigma^h\bigl(\tau(\chi_{\wp})\bigr).\quad \forall h \geq
	1.
\end{eqnarray}
Inserting $h = t f$, one has:
\begin{eqnarray}
\label{4.30}
   \prod_{\wp} \ \tau(\chi_{\wp})^{n^{tf}-1} = 1.
\end{eqnarray}
Let $r | n$ be a prime and $\eu{R} \supset ( r )$ a maximal
ideal. By analogous steps to the proof of Lemma \ref{lcole}, one
eventually shows that:
\begin{eqnarray}
\label{isone}
    \prod_i \ \Bigl( \frac{ \chi_{\wp}( r )}{\eta_i^m } \Bigr)^{ru} \equiv 1
	 \bmod \eu{R}.
\end{eqnarray}
Since $(ru , f ) = 1$, we get $\prod_{\wp} \ \Bigl( \frac{
\chi_{\wp}(r)}{ \eta_i^m } \Bigr) = 1.$ This product of roots of unity
of coprime order can only be $1$ if all factors are $1$ and thus:
\[   \chi_{\wp}(r) = \eta_i^m . \]
The rest of the statement follows by multiplicativity and using
$l_{p(\wp)} (n ) = 1$.
\end{proof}

\section{Deterministic primality test}
The Corollary \ref{jstest}, Theorem \ref{LL} and the certification -
theorem \ref{cert} are used as bases for an explicite primality test,
which proceeds by providing a proof of existence of an \nth{s}
cyclotomic extension of $\id{N}$ for some $s > \sqrt{n}$ such that $t
= \ord_s(n)$ is \textit{small}, de facto
$\id{O}\left(\log(n)^{\log_{(3)}(n)}\right)$.

In all cases, the existence of saturated \nth{p} extensions is
required for all $p | t$. Such an extension or a proof of
compositeness for $n$ can be gained in polynomial time, if one assumes
the existence of some $p$ - power non residues of small height
\cite{BS} -- existence which follows from the GRH. The versions of CPP
based on this assumption are thus probabilistic \textit{Las Vegas}
algorithms; they shall be described with algorithmic details in a
separate paper dealing with implementations. 

The use of GRH is in the case of CPP explicite, in the sense that
the failure to find the required non residues in the expected range
together with an a posteriori proof of primality for $n$, which can be
gained with a variety of methods, would yield a counterexample to the
generalized Riemann hypothesis.

It is however of a certain theoretical interest, that one can prove
also a deterministic version of the Jacobi sum test, one thus that
does not relay upon the existence of saturated extensions. This
version was proposed by Adleman, Pomerance and Rumely in \cite{APR}
and adapted by Lenstra in his exposition \cite{Le1}. Both sources
present the deterministic algorithm as one which is independent of the
\textit{Las Vegas} variant of the Jacobi sum test, and are based on
computation in excessively large extensions. We give here an improved
and simplified version, based on the ideas in \cite{Le1}. Certainly,
the question about the interest for this variant after the AKS test
\cite{AKS} must be addressed. In fact, provided the highly improbable
event occurs, that the Las Vegas version is not sufficient, then the
deterministic version of CPP may still be more efficient than AKS for
larger numbers; this is due both to asymptotic behavior and mostly to
space requirements which are very high for AKS. We add the theory for
the deterministic variant here, for the sake of completeness.

Let thus, as usual, $n$ be an integer to be tested for primality
and $s, t$ integers with $t = \ord_s(n)$ and set

\medskip

$P = \{ \ p | t : p $ is a prime such that $p \ \vert \
\frac{q-1}{\ord_q(n)} $ for some $q | s $ and no saturated \nth{p}
extension of $\id{N}$ is known $ \}.$ 

\medskip

Since the Jacobi sum method can be used for actually constructing
\nth{p} extensions, it follows that in the cases of interest for the
deterministic version, the valuation $v_p(n^{p-1}-1) > 1$ for odd $p$
or $n \equiv 1 \bmod 4$ and $v_2(n^2-1) > 3$ otherwise.

The deterministic test described in \cite{Le1} generalizes the idea of
the Rabin-Miller test. It gives an alternative version of \rf{lpr} in
the $p$-adic numbers $\Z_p$. This leads to proving that the divisors
$r | n$ also lay in some cycles generated by a number $\nu \mod t$,
which can be explicitly constructed: the structure of the criterium is
similar to $(II)$ in Theorem \ref{th1}, replacing $n$ by $\nu$. Since
in general, $\nu \neq n$, the approach paradoxically suggests that no
\nth{s} cyclotomic extensions exist, according to Theorem \ref{th1}.

We consider in depth the case when $p \in P$ is odd or $n \equiv 1
\bmod 4$. The saturation index is in these cases $k_p(n) =
v_p\left(n^{p-1}-1\right)$ and we shall assume that
\[  k_p(n)  = \kappa + 1 > 1 . \]
For such $p$ we let $\id{Q}_p = \left\{ \ \wp = ( \ p^{k(q)}, q ) \ : \ q |
s \ ; \ p^{k(q)} \ \parallel \ (q-1) \ \right\}$ and define $k_m =
\max_{\wp \in \id{Q}_p} \{ k(q) \}$. With this we fix $\zeta =
\xi_{p^{k_m}}$ and for $l < k_m$ we shall assume the compatibility
conditions $\xi_{p^l} = \zeta^{p^{k(q)-l}}$. For a $q | s$ let $\xi =
\xi_q$ be a root of unity, $\Pi \subset G =
\Gal\left(\Q(\xi_q)/\Q\right)$ be the maximal $p$-group, $H = G/\Pi$
and $\eta_q = \sum_{\sigma \in H} \sigma(\xi)$. We shall consider the
rings
\[ \rg{R} = \Z[\zeta] / (n \cdot \Z[\zeta]) \quad \hbox{ and }
 \quad \rg{Q} = \Z[\zeta, \eta_q] / (n \cdot \Z[\zeta, \eta_q]) .\]

Let $\wp = \left(p^{k(q)}, q\right) \in \id{Q}_p$; \ $n^{\varphi\left(
p^{k(q)}\right)} -1 = u(q) \cdot p^{\kappa + k(q)}$ with $ (u(q), p) = 1$,
and fix a character $ \chi = \chi_{\wp} : \ZMs{q} \rightarrow < \zeta
>$.  We assume that
\begin{eqnarray}
\label{minmor}
\left(\tau(\chi)\right)^{n - \sigma_n} = \omega(\chi)^{-n} \in \ < \zeta >
\end{eqnarray}
where $\sigma_a \in \Gal(\Q(\zeta)/\Q) : \zeta \mapsto \zeta^a$, holds
in $\rg{Q}$; a fortiori, if $K = [ \Q(\zeta, \xi) : \Q ]$, we have $
\left(\tau(\chi)\right)^{\frac{n^K-1}{p^{\kappa + k(q)}}} \in \ 
< \zeta >$.  Let $\lambda_q(\chi) = \tau(\chi)^{u(q)} \in \rg{Q}$ and
$\rg{p} =
\begin{cases} p & \hbox{ if $p$ is odd}\\ 4 & \hbox{ otherwise
} \end{cases}.$ With this, we define
\begin{eqnarray}
\label{C.1}
J_{\chi} \subset 1 + \rg{p} \Z_p = \{ \ a \in 1 + \rg{p} \Z_p :
\lambda(\chi)^{a - \sigma_a} \in \ < \zeta > \ \},
\end{eqnarray}

By \rf{minmor}, we have $n^{p-1} \in J_{\chi}$ and thus $J_{\chi}$ is
a non empty subgroup of $\ U_1 = 1 + \rg{p} \Z_p$. The structure of
$U_1$ implies that $J_{\chi} = (1 + p^j)^{\Z_p}$
for a given, yet to determine, positive integer $j$. By analogy  to
the Rabin-Miller test, we let $a_i = 1 + p^i \in \Z_p$, build the
sequence
\[ x_i(q) = (\lambda_q(\chi))^{a_i-\sigma_{a_i}}, \ i = 0, 1, \ldots, k(q)+
\kappa, \] and consider the following conditions:

\begin{itemize}
\item[D1] The halting condition $x_{k(q)+\kappa} \in < \zeta >$ holds. This is
the condition \rf{minmor} and is related to \rf{4.15}.
\item[D2] For $j(q) = \min \{ i : x_i = 1 \}$, $x_{j-1} \in < \zeta >$, if $j >
0$.
\item[D3] If $ \exists \ l \geq 0: x_l \notin < \zeta >$, then for the
maximal such $l$,
\[ \left(x_l - \xi_{p^{k(q)}}^{h}, \ n \right) = 1; \quad 
h = 1, 2, \ldots , p^{k(q)}. \]
\end{itemize}
The first two conditions are Rabin-Miller; by the third, the value of
$j(q)$ in the definition of $J_{\chi}$ is the one determined in D2.
This provides the information that will be used for combining tests.

We now show how this functions. Let $\wp \in \id{Q}_p$, consider a
prime $r | n$ and suppose $r^{p-1} \notin J_{\chi}$. Then, with $j =
j(\wp)$ given by D2, there is an $m \in \Z_p^{\ast}$ with $1 + p^{j-1}
= r^{(p-1)m}$ and the relation \ref{2.38} implies:
\[ x_{j-1} = \tau(\chi)^{r^{(p-1)m} - \sigma_r^{(p-1)m}} =
\chi(r)^{-(p-1)m \cdot r^{(p-1)m}} \mod r \rg{R}.\] This contradicts
condition D3 and thus:
\[   r^{p-1} \in J_{\chi} = (1 + p^j)^{\Z_p}. \]
We shall define $j_p = \max \{ j(\wp) : \wp \in \id{Q}_p \}$ and
choose some $q(p) | s$ such that $q$ gives raise to the maximal value
of $j$, so there is a $\wp = (p^v, q(p)) \in \id{Q}_p$ with $j(\wp) =
j_p$. The condition D3 applied to this particular choice of $\wp$ --
which we shall also refer to as \textit{maximal pair $\wp \in
\id{Q}_p$ } -- implies:
\begin{eqnarray}
\label{mup}
           r^{p-1} = (1 + p^{j_p})^{\mu_p(r)}, \quad \forall \ r \ | \ n,
           \quad \hbox{ with some } \quad \mu_p(r) \in \Z_p.
\end{eqnarray}

\NI Of course, $\Big( n^{p-1} \Big)^{\Z_p} = (1
+p^{1+\kappa})^{\Z_p}$. If $j_p = \kappa + 1$, then $r^{p-1} \in \Big(
n^{p-1} \Big)^{\Z_p}$ for all $r | n$ and consequently, the condition
\rf{lpr} is fulfilled. A saturated \nth{p} extension exists -- albeit,
could not be constructed by the trial and error method of Theorem
\ref{LL}. We deduce from \rf{mup} a condition which is similar to the
one in Lemma \ref{lcole}:
\begin{lemma}
\label{C.2}
Notations being like above, we assume that $n \equiv 1 \bmod 4$ if $p
=2 \in P$. Suppose that the existence of $\mu_p(r)$ in \rf{mup} is
proved by verifying D3 for a maximal $\wp \in \id{Q}_p$ for all $p \in
P$ and that for all $\wp = (p^{k(q)}, q) \in \id{Q}_p$, letting $\chi
= \chi_{\wp} : \ZMs{q} \rightarrow < \zeta_p >$, the condition
\rf{minmor} is verified. Then there is a character:
\begin{eqnarray}
\label{det}
\wh{\chi} : J_{\chi} \rightarrow  < \zeta > \quad \hbox{ with } \quad
\wh{\chi}(r) = \chi(r) \quad \forall \ r | n .
\end{eqnarray}
In particular, $\wh{\chi}(n) = \chi(n)$.
\end{lemma}
\begin{proof}
We may assume that $J_{\chi} = ( a )^{\Z_p}$ with $a = (1 + p^j)$ and
$j \leq j_p$. Let us define $\eta \in < \zeta >$ by the relation
$\eta^{-u a} = \lambda(\chi)^{a - \sigma_a} \in < \zeta > $ and fix
the character $\wh{\chi} : J_{\chi} \rightarrow < \zeta >$ by
$\wh{\chi}(a) = \eta$. If $r | n$ is a prime, by \rf{2.38},
\[ \Big(\tau(\chi)^u \Big)^{r^{p-1}-\sigma_r^{p-1}} = \chi(r)^{-(p-1)u
r^{p-1}} \ \ \bmod r \rg{R},\] while setting $r^{p-1} = (1 +
p^j)^{\mu'} = a^{\mu'}$, with the obvious definition of $\mu'$ in
dependence of $\mu_p(r)$, yields
\[ \Big(\tau(\chi)^u \Big)^{r^{p-1}-\sigma_{r^{p-1}}} =
\lambda(\chi)^{a^{\mu'} - \sigma_{a}^{\mu'}} =
\lambda(\chi)^{(a-\sigma_{a})(\sum_{i=0}^{\mu'-1} \ a^{\mu'-i-1}
\sigma_{a}^i)} = \left(\eta^{-u a}\right)^{\mu' a^{\mu'-1}}.\]
Comparing the last two identities, we find:
\[ \chi(r)^{-(p-1)u r^{p-1}} = \eta^{-\mu' u r^{p-1}} \ \ \mod r \rg{R}.\] 
From $(u r^{p-1}, p) = 1$ and Lemma \ref{lcole} we have $$
\chi(r)^{p-1} = \eta^{\mu'} = \wh{\chi}(a)^{\mu'} =
\wh{\chi}\left(a^{\mu'}\right) = \wh{\chi}\left( r^{p-1} \right).$$
Since $p-1 \in \Z_p^{\ast}$, we also have $\chi(r) = \wh{\chi}(r)$
and, by multiplicativity, $\wh{\chi}(n) = \chi(n)$, which completes
the proof.
\end{proof}
\begin{remark}
It is of practical relevance, to note that all computations can in
fact be performed in the rings $\rg{R} = \Z[\zeta]/\left( n \Z[\zeta]
\right)$, by using multiple Jacobi sums. This is clear for the
verification of \rf{minmor}. In order to determine the value of $j$ in
D2, one has to compute $\left(\tau(\chi)^{a_i-\sigma_{a_i}}\right)^u$
for $a_i = 1+p^i$, and this computation can also be
completed in $\rg{R}$, by definition of the multiple Jacobi sum $J_a(\chi)$.
\end{remark}
\NI Let us introduce the notation $\pi_2(p) = \{ \ q | s \ : \ \exists
\ \wp = \left( \ p^{k(q)}, q \ \right) \in \id{Q}_p \}$ and $\pi_2(P)
= \bigcup_{p \in P} \ \pi_2(p)$. We have the following
deterministic test variant:

\begin{corollary}
Let the notations be like above and suppose that if $2 \in P$ then $n
\equiv 1 \bmod 4$. Suppose that for all $p \not \in P$ and $\wp \in
\id{Q}$ with $\wp = (p^k, q )$, the relation \rf{4.15} holds and that the
existence of the characters $\wh{\chi}$ in Lemma \ref{det} has been
proved for all $\chi = \chi_{\wp}, \wp \in \id{Q}_p$ and $p \in
P$. For all $q | s$, let $ \nu(q)$ be defined by
\begin{eqnarray*}
\chi(\nu(q)) = \begin{cases} \chi(n) & \hbox{ if } q \not \in
\pi_2(P) \\ \wh{\chi_{\wp}}(1+p^{j_p}) & \hbox{ for all } p \in P
\hbox{ with } q \in \pi_2(p), \ \wp = ( p^{j_p}, q ) \in \id{Q}.
\end{cases} 
\end{eqnarray*}
Let $\nu \in \ZMs{s}$ be defined with the Chinese Remainder Theorem,
  by the congruences $\nu \equiv \nu(q) \bmod q$ for all $q | s$.  Then
  all divisors $r | n$ verify $r \in < \nu \bmod s >$.
\end{corollary}
\begin{proof}
Let $r | n$ and $\chi_{\wp}$ be a character, with $\wp = (p^k, q)$; if
$p \not \in P$, then $\chi_{\wp}(\nu) = \chi_{\wp}(n)$ and $\chi(r) =
\chi(\nu)^{l_p(r)}$, as a consequence of Corollary \ref{jstest}. If
$\wp \in \bigcup_{p \in P} \ \id{Q}_p$, then the proof of Lemma
\ref{C.2} implies that $\chi(r) = \chi(\nu)^{\mu_p(r)}$. By choosing
\[ m \equiv \begin{cases}
\mu_p(r) \bmod p^{v_p(t)} &  \hbox{ if } \quad p \in P \\
l_p(r) \bmod p^{v_p(t)} & \hbox{ otherwise },
\end{cases} \]
we find that $\chi(r) = \chi(\nu)^m$ for all characters $\chi \in
\ZMb{s}$. By duality it follows that $r \equiv \nu^m \bmod s$ as
claimed.
\end{proof}

We shall sketch now the case $p = 2$ and $n \equiv 3 \bmod 4$. As
suggested by saturation, we consider here $n^2-1$ instead of $n - 1 =
n^{p-1}-1$ and note that $\Z_p^{\ast} = 3^{\Z_2} \times 5^{\Z_2}$ is
not cyclic any more. For all $q$, one defines like before the
characters $\chi = \chi_{\wp}$ and determines $J_{\chi} \subset
\Z_2^{\ast}$. If $J_{\chi} \neq < n^2 >$, then $n$ is composite, while
for the remaining cases one can define characters $\wh{\chi}$ and show
eventually that an \nth{s} cyclotomic extension of $\id{N}$
exists. There are some technical obstructions \cite{Mi1}, resulting
from the fact that in a first step, only $\wh{\chi}^2$ is naturally
defined and $\wh{\chi}$ having a power of $2$, there is an ambiguity
in its definition. The condition D3 has to be modified and the
ambiguity is removed by considering a $\rho \in \ZMs{s}$ with
${\rho}^2 = 1$ and showing that the possible divisors $r | n$ belong
this time to the set $\{ \nu^k, \rho \nu^k \bmod s : k = 1, 2, \ldots,
t \}$, with $\nu$ defined like in the Corollary. We refer to
\cite{Le1}, \cite{Mi1} for details.

\section{Asymptotics and run times}
In this section we evaluate the asymptotic expected run - time of the
cyclotomy test. We shall use, for ease of notation, the symbol $\id{P}$ for
the set of all rational primes.  The following
theorem is well-known in the context of primality tests \cite{APR},
\cite{CoLe}.
\begin{theorem}[Prachar, Odlyzko, Pomerance]
\label{POP}
There exists an effectively computable positive constant $c$ such that
$\forall \ n > e^e, \ \ \exists t > 0$ satisfying
\begin{eqnarray}
\label{6.2}
t & < & (\log n )^{ c \cdot \log_{(3)}(n)} \quad \mbox{and} \quad
f(t)^2 = \Bigl (\prod_{\{q \in \id{P}, (q-1) \mid t \}} \ q \Bigr)^2 \
> \ n .
\end{eqnarray}
\end{theorem}
Heuristics indicate that the expected value of $c > \log(e)/\log(4)$
and the Theorem shows that one can choose, $(t, s = f(t))$ in the
given range, and then the existence of an \nth{s} cyclotomic extension
can be proved in time polynomial in $t$. The claim follows from $(II)$
of Theorem \ref{th1}. More precisely, if the existence of an \nth{s}
cyclotomic extension is proved by \rf{4.15}, then this relation should
be proved for all pairs $\wp = (p^k, q) \in \id{Q}$, as defined in
Corollary \ref{jstest}. The verification of \rf{4.15} for one fixed
$\wp$ takes $\id{O}^{\sim}\left( p^k \cdot \log(n)^2 \right)$ binary
operations -- with the standard $\id{O}^{\sim}$ notation, in which
factors that are polynomial in $\log(p), \log_{(2)}(n)$ are neglected. We
would wish to deduce some upper bounds on $p^k, q$ and $\sharp \id{Q}$
using the above Theorem. From the prime number Theorem, if $1 < c$ is
such that $\pi(X) < c \cdot \frac{X}{\log(X)}$ for all $X > e^e$, we
have the estimate

\[ \prod_{ p^f < c \cdot \log(X) } \ p^f \ > \ X^{1/2}, \] for all $X
  > e^e$, where $p^f$ are prime powers. Conversely, if $g(X) = \prod_{
    p^f < c \cdot \log(X) } \ p^f$ and $h(Y) = \min \{ X : g(X) >
  Y^{1/2} \}$, the estimate implies:
\begin{eqnarray}
\label{rough}
h(Y) < c \log(Y) \quad \hbox{ and } \quad \pi(h(Y)) < c^2
\frac{\log(Y)}{\log_{(2)}(Y)}, \quad \ \forall  \ Y > 9.
\end{eqnarray}

From this and $q < t$ we deduce that $d(s) < \frac{\log( n
)}{\log_{(2)}(n )}$, where $d(s)$ -- the number of factors of $s =
f(t)$ -- is equal to the number of distinct primes $q$ in the list of
pairs $\id{Q}$. We shall assume here that it is possible to build $t =
\prod_{p^k < B} p^k$ as the product of the first prime powers such
that $f(t) > \sqrt{n}$. This is a hypothesis and not a consequence of
Theorem \ref{POP}. If this holds, it follows from \rf{rough} that for
$\wp = (p^k, q) \in \id{Q}$ we have $p^k < c^2
\log_{(2)}(n)$. Altogether,
\begin{eqnarray}
\label{amalest}
\sharp \id{Q} < c^3 \log(n), \quad p^k < c^2 \cdot \log_{(2)}(n).
\end{eqnarray}
We have the following
\begin{fact}
Let $n, s$ be coprime integers with $n > s > \sqrt{n}$
squarefree. There is a probabilistic Las Vegas algorithm which
requires $\id{O}^{\sim} \left( \log(n)^3 \right)$ binary operations
for {\em proving} the existence of an \nth{s} cyclotomic extension. 
The algorithm generates a certificate for the existence of such extension and the {\em
certificate can be verified}, together with the validity of the Jacobi
sums, in $\id{O}^{\sim} \left( \log(n)^2 \right)$ binary operations.
\end{fact}
\begin{proof}
The proof follows directly from \rf{amalest} and the description of
the algorithm in Section 6. Building up the saturated working
extensions for all primes $p | t$ takes
$\id{O}^{\sim}\left((\log(n)^2\right)$ operations and in the
certificate generation phase, one has to perform an exponentiation with
exponents $O(n)$ in extensions of \textit{small} degree
($O(\log_{(2)}(n))$), for each of $\wp \in \id{Q}$: this leads to the
claimed run time $\id{O}^{\sim}\left(\log(n) \times \log(n)^2
\right)$. The certification requires merely exponents of size
$\id{O}(log_{(2)}(n))$, which explains the verification time, given the
fact that certification of Jacobi sums is negligeable by Lemma \ref{certjac}.
\end{proof}

The operations using superpolynomial time in the CPP primality proofs
are quite elementary: they are the computation of $(2 \cdot \sharp
\id{Q}) \sim \log(n)$ multiple Jacobi sums and the test that $n \bmod
\left( n^k \ \rem \ s \right) \neq 0$ for $k = 2, 3, \ldots, t-1$. Both
operations take $\id{O}^{\sim}( t \log(n) )$ binary operations, and
only the final test is \textit{specific} for $n$; the Jacobi sums can
be reutilized for numerous test and it is conceivable to store large
tables of precomputed sums. Although $t$ and $\log(n)$ are of
different orders of magnitude, we specified the explicite factor
$\log(n)$ for obvious reasons: the exponent of $\log(n)$ in the upper
bound for $t$ diverges so slowly, that it is indicative to know by
what polynomial factor $t$ is multiplied.
\begin{remark}
\label{certtime}
We only estimated the certificates for the existence of \nth{s}
cyclotomic extensions. The existence of such an extension does not
grant primality, and one still has to perform the final trial
divisions \rf{findiv}, requiring a superpolynomial amount of
operations, and for which we did not provide any possible
certification. The interest of CPP certification would be thus rather
theoretical, without a method to \textit{circumvent} \rf{findiv}
completely. 

Such a method is described in \cite{Mi4}, in connection
with {\em dual elliptic primes } and a new algorithm which
intimately combines CPP with ECPP. This combination yields a 
random cubic time primality test with certificates that can be 
verified in quadratic time, being thus the fastest general primality 
test up to date. Like the Atkin version of ECPP, the run time estimates
are based on some heuristics.
\end{remark}

\bibliographystyle{abbrv} 
\bibliography{cpp}

\end{document}